\documentclass[10pt]{article}
\usepackage[utf8]{inputenc}
\usepackage{lmodern}
\usepackage[T1]{fontenc}
\usepackage{amsmath}
\usepackage{amsfonts}
\usepackage{amssymb}
\usepackage{amsthm}
\usepackage{mathrsfs}
\usepackage{stmaryrd}
\usepackage{centernot} 
\usepackage{mathtools}
\usepackage{graphicx}
\usepackage[left=2.9cm,right=2.9cm,top=3cm,bottom=3cm]{geometry}
\usepackage{fancyhdr}
\usepackage{enumitem}
\usepackage{hyperref}
\usepackage{tikz}
\usepackage{dsfont}
\usepackage{multicol}
\usepackage[dvipsnames]{xcolor}
\colorlet{Red}{WildStrawberry}
\usepackage[labelfont=bf,labelsep=colon]{caption}

\usetikzlibrary{calc}

\usepackage[maxbibnames=99]{biblatex}
\renewcommand{\epsilon}{\varepsilon}

\newcommand{\KG}{\mathrm{KG}}

\newcommand{\dd}{\,\mathrm{d} }

\theoremstyle{plain}
\newtheorem{thm}{Theorem}
\newtheorem{lem}[thm]{Lemma}
\newtheorem{prop}[thm]{Proposition}

\theoremstyle{definition}

\theoremstyle{remark}
\newtheorem{rem}[thm]{Remark}

\newcommand\blfootnote[1]{
  \begingroup
  \renewcommand\thefootnote{}\footnote{#1}
  \addtocounter{footnote}{-1}
  \endgroup
}

\begin{document}

\begin{center}
    \large{\textbf{\uppercase{Blow-up, decay, and convergence to equilibrium for focusing damped cubic Klein-Gordon and Duffing equations}}}

    \vspace{\baselineskip}

    \large{\textsc{Thomas Perrin\footnote{\textit{Université de Rennes, ENS Rennes, INRIA, CNRS, IRMAR - UMR 6625, F-35000 Rennes, France.}}}}
\end{center}

\vspace{\baselineskip}

\noindent
\textbf{Abstract.} We study long-time dynamics of the damped focusing cubic Klein-Gordon equation on a compact three-dimensional Riemannian manifold, together with its space-independent reduction, the damped focusing Duffing equation. Under the geometric control condition on the damping and an assumption on the set of stationary solutions, we establish a sharp trichotomy for initial data with energy slightly above that of the ground state: every solution either blows up in finite time, decays exponentially to zero, or converges to a ground state. We provide a complete classification of the Duffing dynamics above the energy of the constant solution, use it to construct Klein-Gordon solutions realising each of the three behaviours in the case of a domain without boundary, and derive a simple spectral criterion ensuring that the ground states are nonconstant--and hence that different types of behaviour can indeed occur for solutions with initial energy above that of the ground state.


\section{Introduction}

The \textit{damped focusing Duffing equation} is given by
\begin{equation}\label{duffing}\tag{$\ast$}
    \left \{
    \begin{array}{rccc}
        \ddot{u} + \gamma \dot{u} + u & = & u^3, & \quad t \geq 0, \\
        \left( u(0), \dot{u}(0) \right) & = & (u_0, u_1), &
    \end{array}
    \right.
\end{equation}
for $(u_0, u_1) \in \mathbb{R}^2$, where $\gamma \geq 0$ is a parameter called the damping. We refer to \eqref{duffing} as \emph{focusing} because some authors (see, for instance, \cite{Duffing_defocusing}) consider the version with the opposite sign on the nonlinearity and refer to that as the damped Duffing equation. 

Our primary motivation for studying \eqref{duffing} stems from its connection to the damped focusing cubic Klein-Gordon equation, as its solutions correspond to space-independent solutions of that equation on a manifold without boundary. More precisely, let $\Omega$ denote the interior of a smooth, compact, connected, three-dimensional Riemannian manifold with (possibly empty) boundary $\partial \Omega$, and let $\Delta$ be the Laplace-Beltrami operator. We denote by $H_0^1(\Omega)$ the closure of $\mathscr{C}_c^\infty(\Omega, \mathbb{R})$ in $H^1(\Omega)$; in particular, $H_0^1(\Omega) = H^1(\Omega)$ if $\partial \Omega = \emptyset$.

The \textit{damped focusing cubic Klein-Gordon equation} is given by
\begin{equation}\label{KG}\tag{$\ast \ast$}
    \left \{
        \begin{array}{rcccl}
            \partial_t^2 u - \Delta u + \gamma \partial_t u + u & = & u^3 & \quad & \text{in } \mathbb{R}_+ \times \Omega, \\
            (u(0), \partial_t u(0)) & = & \left( u_0, u_1 \right) & \quad & \text{in } \Omega, \\
            u & = & 0 & \quad & \text{in } \mathbb{R}_+ \times \partial \Omega,
        \end{array}
    \right.
\end{equation}
for real-valued initial data $\left( u_0, u_1 \right) \in H_0^1(\Omega) \times L^2(\Omega)$ and damping $\gamma \in L^\infty(\Omega, \mathbb{R}_+)$.

\subsection{Definitions and notations for the Duffing equation}

\begin{figure}[ht]
    \centering
    \begin{tikzpicture}[scale=1]
        \def\size{3} 
        \def\exitAbscisse{sqrt(\size*sqrt(2)+1)}
        
        \fill[gray!20]
          plot[domain=-1:1,smooth,variable=\x] 
            ({\x},{(\x*\x -1)/sqrt(2)})
          --
          plot[domain=1:-1,smooth,variable=\x] 
            ({\x},{-(\x*\x -1)/sqrt(2)})
          -- cycle;
        
        \fill[gray!55]
          (\size, -\size) -- (\size, \size) --
          plot[domain=\exitAbscisse:1,smooth,variable=\x]
            ({\x},{(\x*\x -1)/sqrt(2)})
          --
          plot[domain=1:\exitAbscisse,smooth,variable=\x]
            ({\x},{-(\x*\x -1)/sqrt(2)})
          -- cycle;
        
        \fill[gray!55]
          (-\size, \size) -- (-\size, -\size) --
          plot[domain=-\exitAbscisse:-1,smooth,variable=\x]
            ({\x},{-(\x*\x -1)/sqrt(2)})
          --
          plot[domain=-1:-\exitAbscisse,smooth,variable=\x]
            ({\x},{(\x*\x -1)/sqrt(2)})
          -- cycle;
        
        \fill[gray!85]
          plot[domain=\exitAbscisse:1,smooth,variable=\x] 
            ({\x},{(\x*\x -1)/sqrt(2)}) --
          plot[domain=1:-1,smooth,variable=\x] 
            ({\x},{-(\x*\x -1)/sqrt(2)}) --
          plot[domain=-1:-\exitAbscisse,smooth,variable=\x] 
            ({\x},{(\x*\x -1)/sqrt(2)})
          -- cycle;
        
        \fill[gray!85]
          plot[domain=\exitAbscisse:1,smooth,variable=\x] 
            ({\x},{-(\x*\x -1)/sqrt(2)}) --
          plot[domain=1:-1,smooth,variable=\x] 
            ({\x},{(\x*\x -1)/sqrt(2)}) --
          plot[domain=-1:-\exitAbscisse,smooth,variable=\x] 
            ({\x},{-(\x*\x -1)/sqrt(2)})
          -- cycle;
        
        \draw[->] (-\size-0.3,0) -- (\size+0.3,0) node[right] {$u_0$};
        \draw[->] (0,-\size-0.3) -- (0,\size +0.3) node[above] {$u_1$};
        
        \draw[domain=-\exitAbscisse:\exitAbscisse,smooth,variable=\x,black]
          plot ({\x},{(\x*\x -1)/sqrt(2)});
        
        \draw[domain=-\exitAbscisse:\exitAbscisse,smooth,variable=\x,black]
          plot ({\x},{-(\x*\x -1)/sqrt(2)});
        
        \filldraw (1,0) circle (2pt);
        \filldraw (-1,0) circle (2pt);
        
        \node at (1,-0.42) {$+1$};
        \node at (-1,-0.42) {$-1$};
        
        \node[black] at (0.3,0.3) {$\mathcal{K}^+$};
        \node[black] at (2.5,1.5) {$\mathcal{K}^-$};
        \node[black] at (-2.5,1.5) {$\mathcal{K}^-$};
        \node[black] at (0.8,2) {$\mathcal{N}$};
        \node[black] at (0.8,-2) {$\mathcal{N}$};
    \end{tikzpicture}
    \caption{The sets $\mathcal{K}^+$ (light gray), $\mathcal{K}^-$ (medium gray), and $\mathcal{N}$ (dark gray), separated by the curves $\{ E(u_0, u_1) = \frac{1}{4} \}$, given by $u_0 \mapsto \pm \frac{u_0^2 - 1}{\sqrt{2}}$ (in black).}
    \label{fig_K+K-N}
\end{figure}
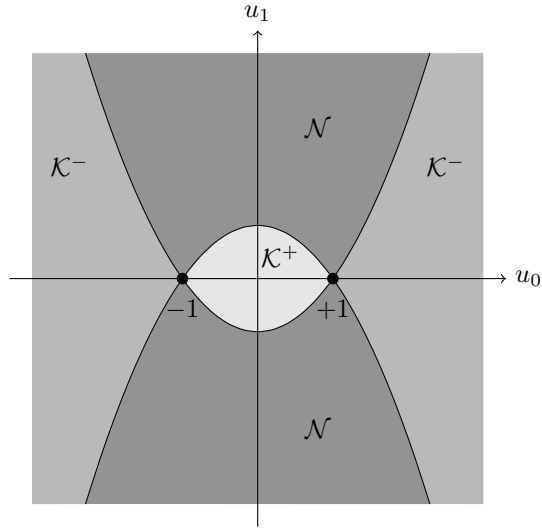

For $(u_0, u_1) \in \mathbb{R}^2$, introduce the energy 
\begin{equation}\label{eq:def_E}
    E(u_0, u_1):= \frac{u_0^2}{2} - \frac{u_0^4}{4} + \frac{u_1^2}{2}.
\end{equation}
The energy of a solution of \eqref{duffing} at time $t \geq 0$ is simply denoted by $E_u(t):= E \left( u(t), \dot{u}(t) \right)$, and one has
\begin{equation}\label{eq_energy_equality}
    E_u(t) = E_u(s) - \gamma \int_s^t \dot{u}(\tau)^2 \dd \tau, 
\end{equation}
for all $t \geq s \geq 0$. Write $\mathbb{R}^2 = \mathcal{K} \sqcup \mathcal{N}$, with
\begin{equation}
    \mathcal{K}:= \left\{ (u_0, u_1) \in \mathbb{R}^2, E(u_0, u_1) < \frac{1}{4} \right\} \ \text{ and } \ \mathcal{N}:= \left\{ (u_0, u_1) \in \mathbb{R}^2, E(u_0, u_1) \geq \frac{1}{4} \right\},
\end{equation}
and further decompose $\mathcal{K} = \mathcal{K}^+ \sqcup \mathcal{K}^-$, with
\begin{equation}
    \mathcal{K}^+:= \left\{ (u_0, u_1) \in \mathcal{K}, \vert u_0 \vert < 1 \right\} \ \text{ and } \ \mathcal{K}^-:= \left\{ (u_0, u_1) \in \mathcal{K}, \vert u_0 \vert > 1 \right\}.
\end{equation}
The sets $\mathcal{K}^+$, $\mathcal{K}^-$ and $\mathcal{N}$ are shown in Figure \ref{fig_K+K-N}.

\subsection{Main results for the Duffing equation}

For solutions with initial data in $\mathcal{K}$, one has the following result, which is the space-independent analogue of the main theorems of \cite{perrin2024damped}, where the damped focusing cubic Klein-Gordon equation is studied. 

\begin{prop}\label{prop_duffing_K}
    Let $(u_0, u_1) \in \mathcal{K}$ and $\gamma \geq 0$.
    \begin{enumerate}[label=(\roman*)]
        \item If $(u_0, u_1) \in \mathcal{K}^-$, then the solution of \eqref{duffing} blows up in finite time. In addition, if $\gamma > 0$, then
        \begin{equation}\label{eq_prop_duffing_K_energy_-infty}
            E_u(t) \xrightarrow{t \rightarrow T^+} - \infty,
        \end{equation}
        where $T^+ < +\infty$ is the maximal time of existence of $u$.
        \item If $(u_0, u_1) \in \mathcal{K}^+$, then the solution of \eqref{duffing} is defined for all $t\geq 0$. In addition, if $\gamma > 0$, then
        \begin{equation}
            \left( u(t), \dot{u}(t) \right) \xrightarrow{t \rightarrow + \infty} (0,0),
        \end{equation}
        at an exponential rate.
    \end{enumerate}
\end{prop}

\begin{figure}[ht]
    \centering
    \begin{tabular}{ccc}
        \includegraphics[width=0.3\textwidth]{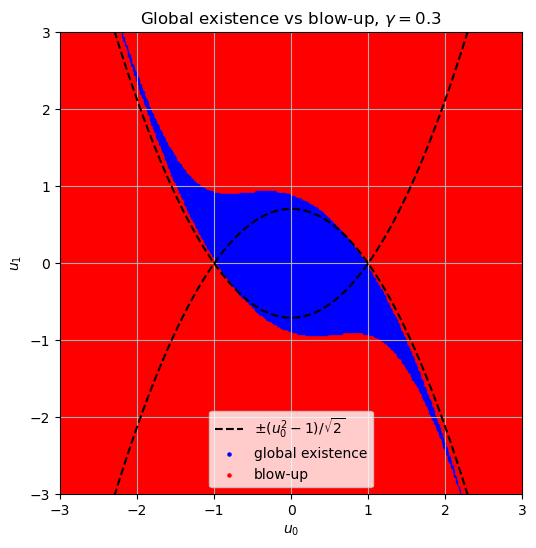} &
        \includegraphics[width=0.3\textwidth]{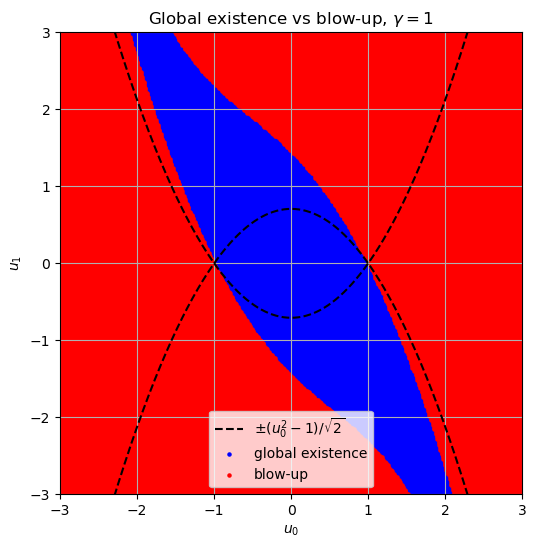} &
        \includegraphics[width=0.3\textwidth]{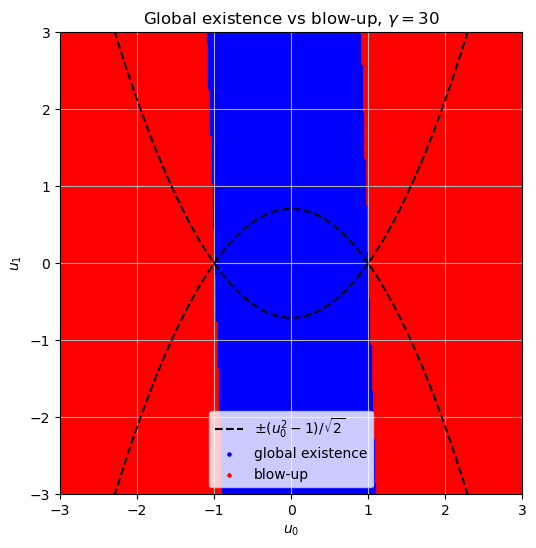}
    \end{tabular}
    \caption{Global existence vs blow-up for three different values of $\gamma$.}
    \label{fig_existence_vs_blowup}
\end{figure}

We provide a proof of Proposition~\ref{prop_duffing_K} below. Before stating our results for solutions with initial data in $\mathcal{N}$, we present some illustrative numerical simulations. In Figure~\ref{fig_existence_vs_blowup}, for three different values of $\gamma$, we solve \eqref{duffing} numerically for many values of $(u_0, u_1) \in \mathbb{R}^2$. A red point is plotted at $(u_0, u_1)$ if the corresponding solution blows up in finite time, and a blue point is plotted if the solution exists globally. Consistently with Proposition~\ref{prop_duffing_K}, we observe that the parameter $\gamma$ has no influence on the blow-up of solutions with initial data in $\mathcal{K}$. However, for initial data in $\mathcal{N}$, the parameter $\gamma$ strongly affects the behaviour of the solution.

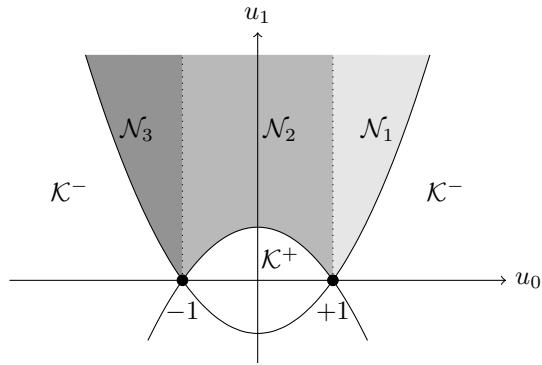
\begin{figure}[ht]
    \centering
    \begin{tikzpicture}[scale=1]
        \def\size{3}
        \def\delta{0.8} 
        \def\exitAbscisse{sqrt(\size*sqrt(2)+1)}
        \def\exitAbscisseLow{sqrt(\delta*sqrt(2)+1)}
        
        \fill[gray!20]
          plot[domain=1:\exitAbscisse,smooth,variable=\x] 
            ({\x},{(\x*\x -1)/sqrt(2)})
          -- (1, \size) -- cycle;
        
        \fill[gray!55]
          plot[domain=-1:1,smooth,variable=\x]
            ({\x},{-(\x*\x -1)/sqrt(2)})
          -- (1,\size) -- (-1, \size) -- cycle;
        
        \fill[gray!85]
          plot[domain=-\exitAbscisse:-1,smooth,variable=\x] 
            ({\x},{(\x*\x -1)/sqrt(2)}) 
          -- (-1, \size) -- cycle;
        
        \draw[->] (-\size-0.3,0) -- (\size+0.3,0) node[right] {$u_0$};
        \draw[->] (0,-\delta-0.3) -- (0,\size +0.3) node[above] {$u_1$};
        
        \draw[domain=-\exitAbscisse:\exitAbscisse,smooth,variable=\x,black]
          plot ({\x},{(\x*\x -1)/sqrt(2)});
        
        \draw[domain=-\exitAbscisseLow:\exitAbscisseLow,smooth,variable=\x,black]
          plot ({\x},{-(\x*\x -1)/sqrt(2)});
        
        \draw[dotted] (1,0) -- (1,\size);
        \draw[dotted] (-1,0) -- (-1,\size);
        
        \filldraw (1,0) circle (2pt);
        \filldraw (-1,0) circle (2pt);
        
        \node at (1,-0.42) {$+1$};
        \node at (-1,-0.42) {$-1$};
        
        \node[black] at (0.3,0.3) {$\mathcal{K}^+$};
        \node[black] at (2.5,1.2) {$\mathcal{K}^-$};
        \node[black] at (-2.5,1.2) {$\mathcal{K}^-$};
        \node[black] at (1.6,2) {$\mathcal{N}_1$};
        \node[black] at (0.3,2) {$\mathcal{N}_2$};
        \node[black] at (-1.6,2) {$\mathcal{N}_3$};
    \end{tikzpicture}
    \caption{The sets $\mathcal{N}_1$ (light gray), $\mathcal{N}_2$ (medium gray), and $\mathcal{N}_3$ (dark gray).}
    \label{fig_K+K-N1N2N3}
\end{figure}

By the transformation $u(t) \longleftrightarrow -u(t)$, it suffices to study (\ref{duffing}) for $(u_0, u_1) \in \mathbb{R} \times [0, +\infty)$. Since the solutions with initial data $(\pm 1, 0)$ are constant, it in fact suffices to consider $(u_0, u_1) \in \mathbb{R} \times (0, +\infty)$. In Figure~\ref{fig_existence_vs_blowup}, we observe that the behaviour of the solution with initial data $(u_0,u_1) \in \mathcal{N}$ depends on the position of $u_0$ relative to $1$ and $-1$. We illustrate separately the phenomena that occur in the three corresponding regions. Write
\begin{equation}
    \mathcal{N} \cap \left\{ (u_0, u_1) \in \mathbb{R}^2, u_1>0 \right\} = \mathcal{N}_1 \sqcup \mathcal{N}_2 \sqcup \mathcal{N}_3,
\end{equation}
with
\begin{equation}
    \mathcal{N}_1:= \left\{ (u_0, u_1) \in \mathcal{N}, u_1>0, u_0 \geq 1 \right\}, \quad 
    \mathcal{N}_2:= \left\{ (u_0, u_1) \in \mathcal{N}, u_1>0, -1 \leq u_0 < 1 \right\},
\end{equation}
and 
\begin{equation}
    \mathcal{N}_3:= \left\{ (u_0, u_1) \in \mathcal{N}, u_1>0, u_0 < -1 \right\}.
\end{equation}
The sets $\mathcal{N}_1$, $\mathcal{N}_2$ and $\mathcal{N}_3$ are shown in Figure \ref{fig_K+K-N1N2N3}.

\begin{figure}[p]
    \centering
    \begin{tabular}{cc}
        \includegraphics[width=0.3\textwidth]{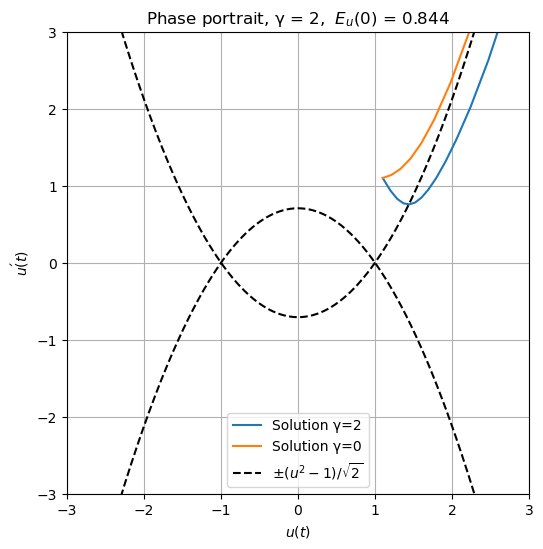} &
        \includegraphics[width=0.3\textwidth]{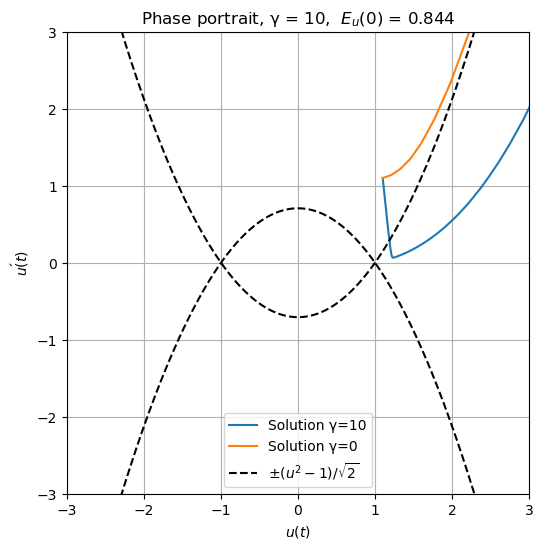} 
    \end{tabular}
    \caption{For $(u_0,u_1) \in \mathcal{N}_1$, the solution blows up for all $\gamma \geq 0$.}
    \label{fig_portait_phase_cas_1}
\end{figure}

First, we consider $(u_0,u_1) \in \mathcal{N}_1$. In Figure \ref{fig_portait_phase_cas_1}, we plot the trajectory of constant energy (corresponding to $\gamma = 0$), and the one with $\gamma > 0$, for two different values of $\gamma$. We observe that if $\gamma > 0$, then the solution exits $\mathcal{N}$, enters $\mathcal{K}^- \cap \{u_0 > 1\}$, and therefore blows up in finite time.

\begin{figure}[p]
    \centering
    \begin{tabular}{ccc}
        \includegraphics[width=0.3\textwidth]{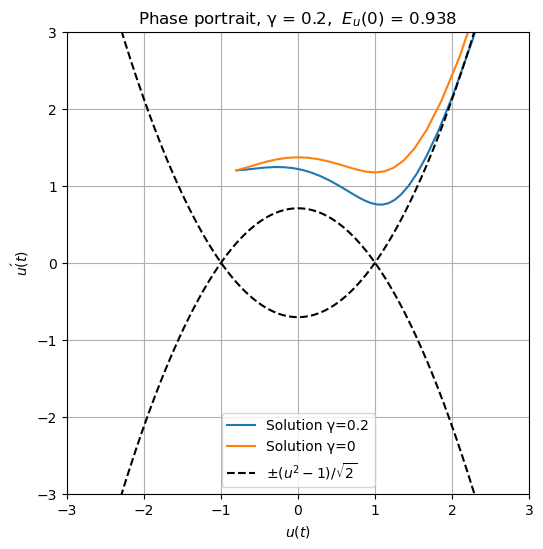} &
        \includegraphics[width=0.3\textwidth]{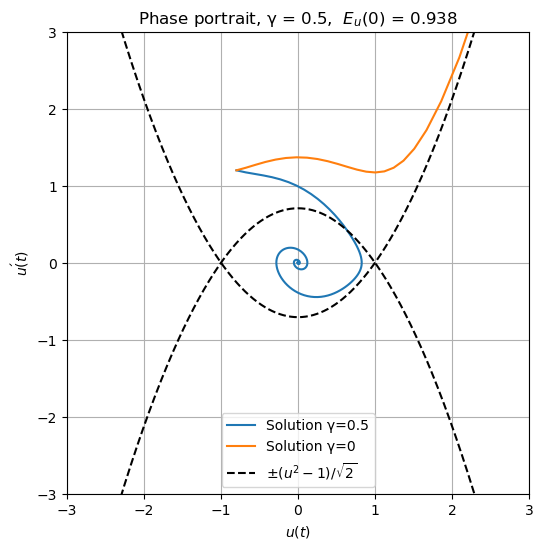} &
        \includegraphics[width=0.3\textwidth]{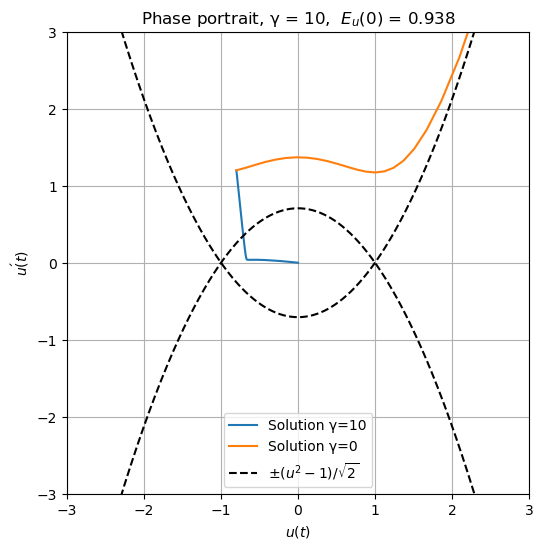} 
    \end{tabular}
    \caption{For $(u_0,u_1) \in \mathcal{N}_2$, the solution blows up for small $\gamma$, and converges to zero for large $\gamma$. One can conjecture that there exists a unique intermediate value $\gamma_0$ such that the corresponding solution converges to $1$.}
    \label{fig_portait_phase_cas_2}
\end{figure}

Second, we consider $(u_0,u_1) \in \mathcal{N}_2$. In Figure \ref{fig_portait_phase_cas_2}, we observe that if $\gamma \geq 0$ is sufficiently small, then the solution enters $\mathcal{N}_1$, and therefore blows up, whereas if $\gamma \geq 0$ is sufficiently large, the solution enters $\mathcal{K}^+$, and therefore converges to zero (at an exponential rate). Assume in addition that $E(u_0, u_1) > \frac{1}{4}$, that is, $(u_0,u_1)$ is in the interior of $\mathcal{N}$, and let $\gamma_0 > 0$ denote the particular value such that for $\gamma < \gamma_0$, the solution blows up, and for $\gamma > \gamma_0$, the solution converges to zero. Then one can conjecture that the solution with damping $\gamma_0$ converges to $(1,0)$.

\begin{figure}[p]
    \centering
    \begin{tabular}{ccc}
        \includegraphics[width=0.3\textwidth]{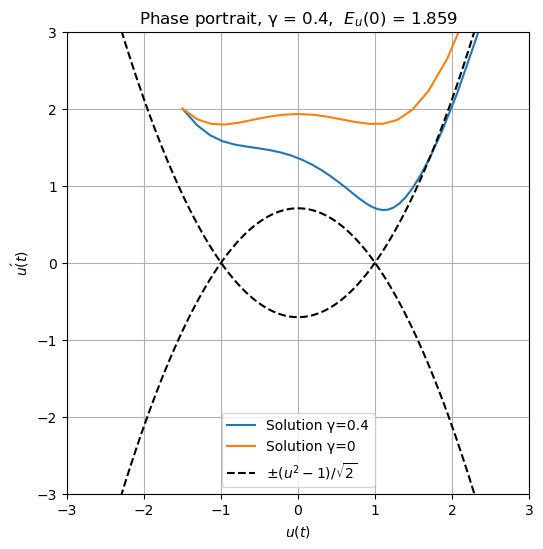} &
        \includegraphics[width=0.3\textwidth]{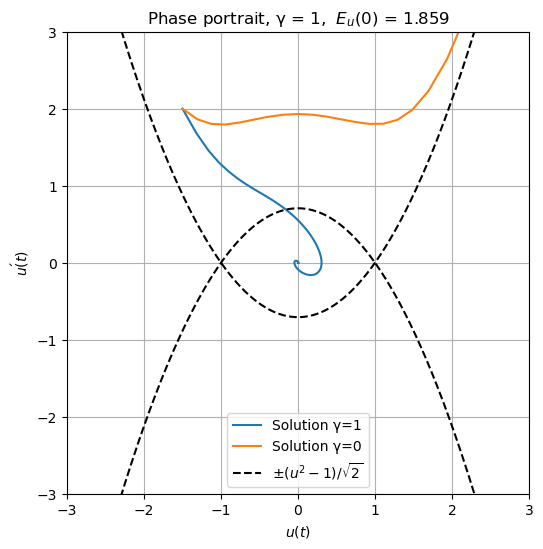} &
        \includegraphics[width=0.3\textwidth]{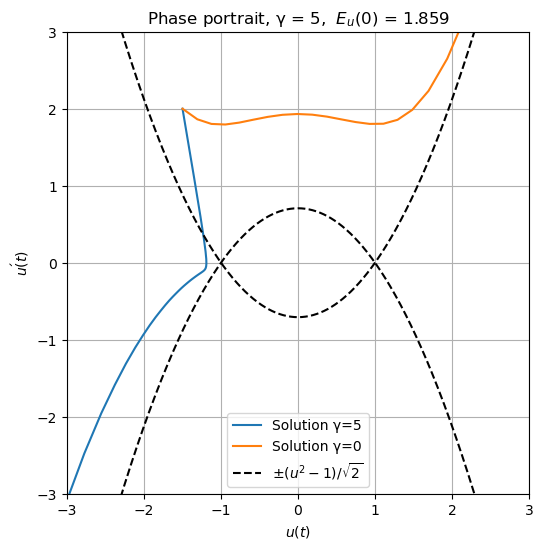} 
    \end{tabular}
    \caption{For $(u_0, u_1) \in \mathcal{N}_3$, the solution blows up for both small and large values of $\gamma$, and converges to zero for some intermediate values of $\gamma$. One can conjecture that there exist two intermediate values $\gamma_1 > \gamma_0$ such that the corresponding solutions converge to $-1$ and $1$ respectively.}
    \label{fig_portait_phase_cas_3}
\end{figure}

Finally, we consider $(u_0, u_1) \in \mathcal{N}_3$. In Figure~\ref{fig_portait_phase_cas_3}, we observe that if $\gamma \geq 0$ is sufficiently small, then the solution enters $\mathcal{N}_2$, and therefore either blows up (if $\gamma$ is small enough) or converges to zero (if $\gamma$ is slightly larger). On the other hand, if $\gamma$ is large, the solution enters $\mathcal{K}^-$ and thus blows up in finite time. As above, if $E(u_0, u_1) > \frac{1}{4}$ and if $\gamma_0 < \gamma_1$ denote the two particular values of the damping, one can conjecture that the solution with damping $\gamma_0$ converges to $(1,0)$, while the solution with damping $\gamma_1$ converges to $(-1,0)$.

These observations are made precise in the following theorem, which is our main result about the damped Duffing equation.

\begin{thm}\label{thm_duffing_N}
    Let $(u_0, u_1) \in \mathcal{N}$, with $u_1>0$, and let $\gamma \geq 0$.
    \begin{enumerate}[label=(\roman*)]
    
        \item If $(u_0, u_1) \in \mathcal{N}_1$, then the solution of \eqref{duffing} blows up in finite time. In addition, if $\gamma > 0$, then the solution enters $\mathcal{K}^-\cap \{ u_0 > 1 \}$.
        
        \item Assume $(u_0, u_1) \in \mathcal{N}_2$. Then, there exists $\gamma_0 \geq 0$ satisfying the following property: first, if $\gamma < \gamma_0$, then the solution of \eqref{duffing} enters $\mathcal{N}_1$ and thus blows up in finite time, second, if $\gamma = \gamma_0$, then 
        \begin{equation}
            \left( u(t), \dot{u}(t) \right) \xrightarrow{t \rightarrow + \infty} (1,0),
        \end{equation}
        at an exponential rate, and third, if $\gamma > \gamma_0$, then 
        \begin{equation}
            \left( u(t), \dot{u}(t) \right) \xrightarrow{t \rightarrow + \infty} (0,0),
        \end{equation}
        at an exponential rate. If in addition $E(u_0, u_1) = \frac{1}{4}$, then $\gamma_0 = 0$.
        
        \item Assume $(u_0, u_1) \in \mathcal{N}_3$, with $E(u_0, u_1) > \frac{1}{4}$. Then, there exist $\gamma_1 > \gamma_0 \geq 0$ satisfying the following property: first, if $\gamma < \gamma_0$, then the solution of \eqref{duffing} enters $\mathcal{N}_1$ and thus blows up in finite time, second, if $\gamma = \gamma_0$, then 
        \begin{equation}
            \left( u(t), \dot{u}(t) \right) \xrightarrow{t \rightarrow + \infty} (1,0),
        \end{equation}
        at an exponential rate, third, if $\gamma_1 > \gamma > \gamma_0$, then 
        \begin{equation}
            \left( u(t), \dot{u}(t) \right) \xrightarrow{t \rightarrow + \infty} (0,0),
        \end{equation}
        at an exponential rate, fourth, if $\gamma = \gamma_1$, then 
        \begin{equation}
            \left( u(t), \dot{u}(t) \right) \xrightarrow{t \rightarrow + \infty} (-1,0),
        \end{equation}
        at an exponential rate, and fifth, if $\gamma > \gamma_1$, then the solution enters $\mathcal{K}^-\cap \{ u_0 < -1 \}$, and therefore blows up in finite time.

        \item Assume $(u_0, u_1) \in \mathcal{N}_3$, with $E(u_0, u_1) = \frac{1}{4}$. If $\gamma=0$ then
        \begin{equation}
            \left( u(t), \dot{u}(t) \right) \xrightarrow{t \rightarrow + \infty} (-1,0),
        \end{equation}
        at an exponential rate, and if $\gamma > 0$, then the solution enters $\mathcal{K}^-\cap \{ u_0 < -1 \}$, and therefore blows up in finite time.
    \end{enumerate}
\end{thm}

\subsection{Definitions and notations for the Klein-Gordon equation}

We write $\left\Vert u_0 \right\Vert_{H_0^1}^2 = \int_\Omega \left(\left\vert \nabla u_0 \right\vert^2 + \vert u_0 \vert^2 \right) \dd x$, and we define two functionals $J_\KG$ and $E_\KG$, which are respectively called \emph{static energy} and \emph{energy}, by
\begin{equation}
    J_\KG(u_0) = \frac{1}{2} \left\Vert u_0 \right\Vert_{H_0^1}^2 - \frac{1}{4} \left\Vert u_0 \right\Vert_{L^4}^4 \quad \text{ and } \quad E_\KG\left( u_0, u_1 \right) = J_\KG(u_0) + \frac{1}{2} \left\Vert u_1 \right\Vert_{L^2}^2.
\end{equation}
The energy of a solution of (\ref{KG}) at time $t \in \mathbb{R}$ is defined by $E_\KG\left( u(t), \partial_t u(t) \right)$, and satisfies
\begin{equation}
    E_\KG\left( u(t), \partial_t u(t) \right) = E_\KG\left( u_0, u_1 \right) - \int_0^t \int_\Omega \gamma(x) \left\vert \partial_t u(s, x) \right\vert^2 \dd x \dd s.
\end{equation}
In the case $\partial \Omega = \emptyset$, note that if $\left( u_0, u_1 \right) \in H^1(\Omega) \times L^2(\Omega)$ is constant, that is, if $\left( u_0(x), u_1(x) \right) = \left( u_0, u_1 \right) \in \mathbb{R}^2$ for all $x \in \Omega$, then $E_\KG\left( u_0, u_1 \right) = \vert \Omega \vert E\left( u_0, u_1 \right)$, where $E$ is defined by \eqref{eq:def_E} and $\vert \Omega \vert$ is the volume of $\Omega$.

Let $Q$ be a ground state of (\ref{KG}), that is, a stationary solution of (\ref{KG}) of minimal energy 
\begin{equation}
    E_\KG(Q, 0) = J_\KG(Q) =: d > 0.
\end{equation}
See Lemma \ref{lem_ground_state_properties} for a precise definition. For $u_0 \in H_0^1(\Omega)$, write 
\begin{equation}\label{eq_def_K}
    K(u_0) = \left\Vert u_0 \right\Vert_{H_0^1}^2 - \left\Vert u_0 \right\Vert_{L^4}^4,
\end{equation}
and set
\begin{equation}
    \left \{
        \begin{array}{c}
            \mathcal{K}_\KG^+ = \left\{ \left( u_0, u_1 \right) \in H_0^1(\Omega) \times L^2(\Omega), E_\KG\left(u_0, u_1 \right) < d, K(u_0) \geq 0 \right\}, \\
            \mathcal{K}_\KG^- = \left\{ \left( u_0, u_1 \right) \in H_0^1(\Omega) \times L^2(\Omega), E_\KG\left(u_0, u_1 \right) < d, K(u_0) < 0 \right\}.
        \end{array}
    \right.
\end{equation}

For $\omega \subset \Omega$, we say that $\omega$ satisfies the \emph{Geometric Control Condition} (in short, GCC) if there exists $L > 0$ such that any generalised geodesic of $\Omega$ of length $L$ meets the set $\omega$. For the definition of generalised geodesics, we refer to \cite{MelroseSjostrand}. We always assume that no generalised geodesic has a contact of infinite order with $\partial \Omega$ (see \cite{BLR} for some details about this assumption). By Theorems 2 and 4 of \cite{perrin2024damped}, the sets $\mathcal{K}_\KG^+$ and $\mathcal{K}_\KG^-$ are stable under the forward flow of \eqref{KG}; solutions initiated in $\mathcal{K}_\KG^-$ blow up in finite time, while those initiated in $\mathcal{K}_\KG^+$ exist globally and, moreover, are stabilised when $\gamma$ satisfies the GCC. 

\subsection{Main results for the Klein-Gordon equation}

The following theorem is a consequence of the results established in \cite{perrin2024damped}, and is in fact valid if the first line of \eqref{KG} is replaced by $\partial_t^2 u - \Delta u + \gamma \partial_t u + \beta u = u^3$, for any $\beta \in \mathbb{R}$ such that the Poincaré's inequality
\begin{equation}
   \int_\Omega \vert u_0 \vert^2 \dd x \leq C \int_\Omega \left(\left\vert \nabla u_0 \right\vert^2 + \beta \vert u_0 \vert^2 \right) \dd x
\end{equation}
holds for all $u \in H_0^1(\Omega)$, for some $C > 0$, and if the definition of $\Vert u \Vert_{H_0^1}$ is modified accordingly (see \cite{perrin2024damped}).

\begin{thm}\label{thm:conjecture_ok_strong_assumptions}
    Suppose that $\gamma(x) \geq \alpha$ for almost every $x \in \omega$, where $\alpha > 0$ is a constant and $\omega \subset \Omega$ is an open set satisfying the GCC. In addition, assume that the set of ground states is at most countable, and that there exists $\epsilon > 0$ such that no stationary solution $w$ of \eqref{KG} has energy $J(w) \in (d, d+\epsilon]$. Then, for all $\left( u_0, u_1 \right) \in H_0^1(\Omega) \times L^2(\Omega)$ with $E_\KG(u_0, u_1) \leq d + \epsilon$, one of the following three alternatives holds.
    \begin{enumerate}[label=(\roman*)]
        \item The solution $u$ blows up in finite time.
        \item The solution $u$ is stabilised, that is,
        \begin{equation}
            \left( u(t), \partial_t u(t) \right) \xrightarrow{t \rightarrow + \infty} (0, 0) \quad \text{ in } H_0^1(\Omega) \times L^2(\Omega),
        \end{equation}
        at an exponential rate.
        \item There exists a ground state $Q$ of \eqref{KG} such that
        \begin{equation}
            \left( u(t), \partial_t u(t) \right) \xrightarrow{t \rightarrow + \infty} (Q, 0) \quad \text{ in } H_0^1(\Omega) \times L^2(\Omega).
        \end{equation}
    \end{enumerate}
\end{thm}

\begin{proof}
    Let $\left( u_0, u_1 \right) \in H_0^1(\Omega) \times L^2(\Omega)$ satisfy $E_\KG(u_0, u_1) \leq d + \epsilon$. If $u$ blows up in finite time, then \textit{(i)} holds, so we may assume that $u$ is defined on $[0, + \infty)$. If there exists $t \geq 0$ such that $E_\KG\left( u(t), \partial_t u(t) \right) < d$, then since $u$ does not blow up, Theorem 2 of \cite{perrin2024damped} implies that $u$ enters $\mathcal{K}_\KG^+$, and hence \textit{(ii)} holds by Theorem 4 of \cite{perrin2024damped}. Finally, assume that $E_\KG\left( u(t), \partial_t u(t) \right) \in [d, d+\epsilon)$ for all $t\geq 0$. Set 
    \begin{equation}
        E_\infty = \lim_{t \rightarrow + \infty} E_\KG\left( u(t), \partial_t u(t) \right) \in [d, d+\epsilon].
    \end{equation}
    Since there is no stationary solution of \eqref{KG} with energy in $(d, d+\epsilon]$, and since the set of ground states is at most countable, Theorem 5 of \cite{perrin2024damped} implies that \textit{(iii)} holds.
\end{proof}

Theorem \ref{thm_duffing_N} and the methods developed in its proof can be used to construct solutions of \eqref{KG} with interesting behaviours. By definition, $J_\KG(1) = J_\KG(x \mapsto 1) = \frac{\vert \Omega \vert}{4}$.

\begin{thm}\label{thm:construction_solution_particuliere}
    Assume that $\partial \Omega = \emptyset$ and that $\gamma$ is a positive constant. Then, there exist two nonempty open subsets $\mathcal{O}_1$, $\mathcal{O}_2$ of 
    \begin{equation}
        \left\{ \left( u_0, u_1 \right) \in H_0^1(\Omega) \times L^2(\Omega), E_\KG(u_0, u_1) > J_\KG(1) \right\}
    \end{equation}
    such that for all $\left( u_0, u_1 \right) \in \mathcal{O}_1$, \textit{(i)} holds, and for all $\left( u_0, u_1 \right) \in \mathcal{O}_2$, \textit{(ii)} holds. In addition, there exists $\left( u_0, u_1 \right) \in H^1(\Omega) \times L^2(\Omega)$ with $E_\KG(u_0, u_1) > J_\KG(1)$ such that the solution $u$ satisfies
    \begin{equation}
        \left( u(t), \partial_t u(t) \right) \xrightarrow{t \rightarrow + \infty} (+1, 0) \quad \text{ in } H^1(\Omega) \times L^2(\Omega),
    \end{equation}
    at an exponential rate. Likewise, there exists $(u_0, u_1)$, with $E_\KG(u_0, u_1) > J_\KG(1)$, for which the solution converges exponentially to $(-1, 0)$.
\end{thm}

Note that Theorem \ref{thm:construction_solution_particuliere} is true without any assumption on the ground states. If $\partial \Omega \neq \emptyset$, then the proof of Theorem \ref{thm:construction_solution_particuliere} can be adapted away from the boundary if $\gamma \in L^\infty(\Omega, \mathbb{R}_+)$ is constant on an open subset of $\Omega$, using the finite speed of propagation. However, this case is of limited interest, since the condition $E_\KG(u_0, u_1) > J_\KG(1)$ would no longer be satisfied.

If $\Omega$ is such that the ground states are the constant functions $\pm 1$, then Theorem \ref{thm:construction_solution_particuliere} provides situations in which all three possibilities of Theorem \ref{thm:conjecture_ok_strong_assumptions} may occur. If $d = J_\KG(Q) \leq J_\KG(1)$ for all ground states $Q$, then the solutions given by Theorem \ref{thm:construction_solution_particuliere} have an initial energy above that of the ground state. We give a sufficient condition for the condition $d = J_\KG(Q) \leq J_\KG(1)$ to hold.

\begin{prop}\label{prop:Q_nonconstant_lambda_1_smaller_than_2}
    Assume that $\partial \Omega = \emptyset$, and that $\lambda_1 < 2$, where $\lambda_1$ is the first nonzero eigenvalue of $-\Delta$. Then for every ground state $Q$ of \eqref{KG}, one has $J_{\KG}(Q) < J_{\KG}(1)$, that is, the constant solutions $\pm 1$ are stationary solutions of \eqref{KG} with energy strictly above the ground state energy. 
\end{prop}

The question of determining whether ground states are $\pm 1$ if $\lambda_1 \geq 2$ seems to depend strongly on the global geometry of $\Omega$. See below for a few bibliographic remarks about this subject.

\subsection{Bibliographic comments}

The study of the long-time behaviour of solutions and of stationary solutions of nonlinear wave equations has received considerable attention.
We provide here a non-exhaustive list of some references about these topics.

\paragraph{Long-time behaviour of solutions.}
The explicit partition $\mathcal{K}_\KG^+ \sqcup \mathcal{K}_\KG^-$ of the set of initial data below the ground-state energy into blow-up and global solutions has been known since the work of Payne and Sattinger \cite{Payne-Sattinger}. The article \cite{perrin2024damped} by the present author extends the results of Payne and Sattinger to the case of the damped equation, and proves the stabilisation of global solutions below the ground-state energy, under the geometric control condition, in the spirit of \cite{Joly-Laurent}. An analogous definition of $\mathcal{K}_\KG^\pm$ can also be made for a nonlinearity of the form $u^p$, and for the wave equation instead of the Klein-Gordon equation. A major breakthrough was achieved in \cite{KenigMerle}, in which the authors proved that for the critical wave equation on $\mathbb{R}^n$, an explicit dichotomy between blow-up and scattering solutions occurs. The corresponding result for the Klein-Gordon equation was later established in \cite{IbrahimMasmoudiNakanishi}. 

The behaviour of solutions with energy exactly equal to that of the ground state was then studied in \cite{DuyckaertsMerle}, for the critical wave equation on $\mathbb{R}^N$, where the authors proved the existence of two heteroclinic solutions: both scatter to the ground state in positive time, and in negative time respectively scatter to zero or blow up. The authors of \cite{NakanishiSchlagArticleRAD} and \cite{NakanishiSchlagArticleNONRAD} then studied solutions of the cubic Klein-Gordon equation on $\mathbb{R}^3$ with energy slightly above that of the ground state. They proved the existence of the same heteroclinic solutions as in \cite{DuyckaertsMerle}, as well as a new type of behaviour, namely blow-up in positive time and scattering to zero in negative time. The analogous results for the wave equation were later obtained in \cite{krieger_global_2012, krieger_global_2013}. Related properties can be found in \cite{krieger_threshold_2014, krieger_center-stable_2013, Duyckaerts2024} and the references therein. For a pedagogical presentation of these results, we refer to the book \cite{NakanishiSchlagBOOK}. The one-dimensional case was studied in \cite{krieger_global_2012-1D}.

Damped nonlinear wave equations have also been studied by several authors. In \cite{BurqRaugelSchlag}, in the case of the subcritical Klein-Gordon equation on $\mathbb{R}^d$ with constant damping, the authors proved that radial solutions either blow up or converge to an equilibrium point. The authors of \cite{LiZaho2020} showed that solutions either blow up, exist globally and are unbounded, or converge to a sum of equilibrium points--still in the constant damping case, but without the radial assumption. On a bounded domain, blow-up of certain solutions of a (possibly strongly) damped equation under the ground-state energy was proved in \cite{pucci_global_1998}, and some extensions allowing the precise value of the ground-state energy can be found in \cite{vitillaro_global_1999}. Further results for both strongly and weakly damped equations were given in \cite{gazzola_global_2006}. Finally, the damped defocusing and the undamped focusing Duffing equations were studied in \cite{Duffing_defocusing} and \cite{Duyckaerts2024}, respectively.

\paragraph{Ground states and related properties.}
Note that the results vary significantly depending on the geometry of $\Omega$ and on the precise form of the equation (for instance, on whether one studies $-\Delta u = u^p$ or $-\Delta u + u = u^p$). A proof of the existence of ground states was given, for instance, in \cite{Payne-Sattinger}. In many situations, ground states do not change sign (see, for example, \cite[Theorem 2.3]{Payne-Sattinger}). However, the uniqueness result stated in \cite{Payne-Sattinger} is false. Additional existence results for stationary solutions in $\mathbb{R}^n$ were proved in \cite{berestycki_methode_1980} using a different method.

In \cite{Coffman1972}, the author proved the uniqueness of the positive radial stationary solution of \eqref{KG} in the case $\Omega = \mathbb{R}^3$. Several generalisations were obtained in \cite{chen_uniqueness_1991}. In \cite{Ni1983, Bandle1989}, the authors established uniqueness results for radial solutions on a ball or an annulus. For uniqueness results for radial solutions on possibly unbounded radially symmetric domains, in the case of Dirichlet-Neumann boundary conditions, see \cite{Kwong1989}. Later, the author of \cite{tang_uniqueness_2003} provided a simpler proof of the uniqueness of radial solutions on (possibly unbounded) annuli.

In \cite{Gidas1979}, for a nonlinear wave equation on a ball, the authors proved the uniqueness of the positive stationary solution, and showed that it is radial (see equation (2.8) in \cite{Gidas1979}). Several symmetry properties for stationary solutions on bounded domains were established in \cite{li_monotonicity_1991}. In the famous article \cite{brezis_positive_1983}, the authors investigated the problem $-\Delta u = u^p + f(x,u)$ on a bounded domain, and in particular the model $-\Delta u = u^p + \lambda u$, when $p$ is the critical Sobolev exponent. On page 453, they proved the existence of nonradial ground states of $-\Delta u = u^q$ on an annulus when $q$ is close to the critical Sobolev exponent. In particular, the set of ground states may be uncountable. The author of \cite{lin_existence_1992} studied the existence of nonradial stationary solutions on an annulus, depending on its thickness. In \cite{lin_asymptotic_1995}, the same author studied the case of an annulus of fixed thickness and large radius. A uniqueness result for wave equations with subcritical nonlinearity on bounded domains was proved in \cite{zou_effect_1994}, under the assumption that the domain is ``boundedly different from a ball'' in a specific sense. The works \cite{byeon_existence_1997, byeon_effect_2000} (with the addendum \cite{byeon_addendum_2001}) studied the number of positive solutions on an annulus with subcritical nonlinearity. The author proved that the dimension of the space has a strong impact on the results and showed that the case $n=3$ is surprisingly different from the others. The influence of symmetries on the solutions was also investigated.

On $\mathbb{R}^n$ or on annular domains, the existence and uniqueness of radial stationary solutions with a prescribed number of zeros were studied in \cite{jones_infinitely_1986, seda_properties_1990, grillakis_existence_1990, mcleod_radial_1990, troy_existence_2005, yarur_uniqueness_2009, tanaka_uniqueness_2016}. We note that on non-compact domains, situations may arise in which ground states do not exist, as demonstrated in \cite{Esteban1982} and \cite{Benci1987}. Some existence results for unbounded domains were provided in \cite{noussair_existence_1979, bahri_existence_1997}. Radial stationary solutions on $\mathbb{R}^n$ that are singular at the origin were constructed in \cite{johnson_singular_1994}.

\subsection{Outline of the article} 

In Section \ref{sec:duffing}, we prove our main results about the Duffing equation, Proposition \ref{prop_duffing_K} and Theorem \ref{thm_duffing_N}, together with an additional result (see Lemma \ref{lem:gamma_fix_u_1_varies}), which will be used to prove Theorem \ref{thm:construction_solution_particuliere}. Section~\ref{sec:KG} begins with the well-posedness of solutions of the Klein-Gordon equation and a source-to-solution continuity result. We then recall the definition and some properties of ground states and establish Proposition~\ref{prop:Q_nonconstant_lambda_1_smaller_than_2} about the ground states and the constant solutions $\pm 1$. Finally, we prove Theorem~\ref{thm:construction_solution_particuliere}.

\subsection{Notation} 

When an estimate of the form $A(x) \leq B(x)$ is valid for all $x$ in a set $E$ (or, more generally, when a property depending on certain variables holds throughout a given range), we may use the shorthand notation
\begin{equation}
    A(x) \leq B(x), \quad x \in E.
\end{equation}
When the maximal time of existence $T^+$ of a solution is finite, we often write $A(t) \xrightarrow{t \rightarrow T^+} B$ as a shorthand for
\begin{equation}
    A(t) \xrightarrow[t < T^+]{t \rightarrow T^+} B.
\end{equation}

\section{Behaviour of solutions of the Duffing equation}\label{sec:duffing}

\subsection{Solutions of the Duffing equation of low energy}

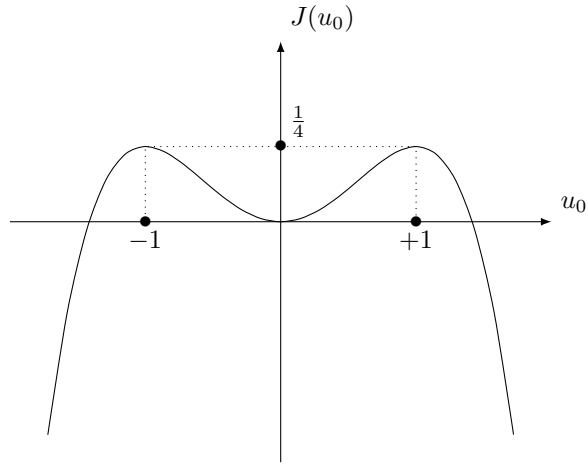
\begin{figure}[ht]
    \centering
    \begin{tikzpicture}[yscale=4,xscale=1.8,>=latex]
        \draw[->] (-2,0) -- (2,0) node[above right] {$u_0$};
        \draw[->] (0,-0.8) -- (0,0.6) node[above right] {$J(u_0)$};
        
        \draw[domain=-1.72:1.72,smooth,variable=\x,black]
          plot ({\x},{\x*\x / 2 - \x*\x*\x*\x / 4});
        
        \draw[dotted] (1,0) -- (1,0.25) -- (-1,0.25) -- (-1,0);
        \node[above right] at (0,0.25) {$\frac{1}{4}$};
        
        \node at (0,0.25) {$\bullet$};
        \node at (1,0) {$\bullet$};
        \node at (-1,0) {$\bullet$};
        
        \node[below] at (1,0) {$+1$};
        \node[below] at (-1,0) {$-1$};
    \end{tikzpicture}
    \caption{Graph of the function $J$.}
    \label{fig_graph_J}
\end{figure}

Here, we prove Proposition \ref{prop_duffing_K}. Set $J(u_0) = \frac{(u_0)^2}{2} - \frac{(u_0)^4}{4}$, so that $E(u_0,u_1) = J(u_0) + \frac{u_1^2}{2}$. The graph of $J$ is shown in Figure \ref{fig_graph_J}.

\paragraph{Solutions with initial data in $\mathcal{K}^-$.} Here, we prove Proposition \ref{prop_duffing_K} \emph{(i)}. Let $(u_0, u_1) \in \mathcal{K}^-$, $\gamma \geq 0$, and write $T^+ \in (0, +\infty]$ for the maximal time of existence of the solution $u$ of \eqref{duffing}. Note that the function $t \mapsto J(u(t))$ is continuous, and satisfies 
\begin{equation}
    J(u(t)) \leq E_u(t) \leq E(u_0, u_1) < \frac{1}{4}.
\end{equation}
Hence, as $\vert u(0) \vert  > 1$, there exists $\delta > 0$ such that 
\begin{equation}\label{eq:proof:potential_well}
    u(t)^2 \geq 1 + \delta, \quad t \in [0, T^+).
\end{equation}
Hence, the following lemma implies that $u$ blows up in finite time.

\begin{lem}\label{lem_blow-up_geq_1_plus_delta}
    Let $(u_0, u_1) \in \mathbb{R}^2$, $\gamma \geq 0$, and write $T^+ \in (0, +\infty]$ for the maximal time of existence of the solution $u$ of \eqref{duffing}. If there exists $\delta > 0$ such that 
    \begin{equation}\label{eq_proof_blow_up_K-_1}
        u(t)^2 \geq 1 + \delta, \quad t \in [0, T^+),
    \end{equation}
    then $T^+ < + \infty$, that is, $u$ blows up in finite time.
\end{lem}

\begin{proof}\label{lem_u_>_1_blowup}
    Assume by contradiction that $T^+ = + \infty$. For $t \geq 0$, set $M(t) = u(t)^2$. One has $M^\prime = 2 \dot{u} u$ and 
    \begin{equation}\label{eq_M_prime_prime}
        M^{\prime \prime} = 2 \left( u^4 - u^2 - \gamma u \dot{u} + \dot{u}^2 \right).
    \end{equation}

    First, we prove 
    \begin{equation}\label{eq_proof_blow_up_K-_2}
        M(t) \xrightarrow{t \rightarrow +\infty} +\infty.
    \end{equation}
    Using (\ref{eq_proof_blow_up_K-_1}) and (\ref{eq_M_prime_prime}), one finds
    \begin{equation}
        M^{\prime \prime} + \gamma M^\prime \geq 2 \delta u^2 \geq 2 \delta (1 + \delta) > 0,
    \end{equation}
    implying that there exists $T > 0$ such that
    \begin{equation}
        M^{\prime}(t) + \gamma M(t) \geq c t, \quad t \geq T,
    \end{equation}
    for some $c > 0$. Integrating, one obtains
    \begin{equation}
        M(t) \geq M(T) e^{\gamma(T - t)} + c e^{- \gamma t} \int_{T}^t \tau e^{\gamma\tau} \dd \tau, \quad t\geq T,
    \end{equation}
    and this gives (\ref{eq_proof_blow_up_K-_2}).

    Second, we obtain a contradiction.
    For $\epsilon > 0$, one has
    \begin{equation}
        M^{\prime \prime} \geq 2 u^4 - 2 u^2 - \frac{C_0}{\epsilon} u^2 - \epsilon \dot{u}^2 + 2 \dot{u}^2,
    \end{equation}
    for some $C_0 > 0$. 
    Set (for instance) $\epsilon = \frac{1}{2}$. By (\ref{eq_proof_blow_up_K-_2}), there exists $T_1 \geq 0$ such that 
    \begin{equation}
        M^{\prime \prime}(t) \geq \frac{3}{2} u(t)^4 + \frac{3}{2} \dot{u}(t)^2 = \frac{9}{2} \dot{u}(t)^2 + 3 u(t)^2 - 6E_u(t), \quad t \geq T_1.
    \end{equation}
    As the energy is nonincreasing, one obtains
    \begin{equation}
        M^{\prime \prime}(t) \geq \frac{9}{2} \dot{u}(t)^2 + 3 u(t)^2 - 6 E_u(0) \geq \frac{9}{2} \dot{u}(t)^2, \quad t \geq T_2,
    \end{equation}
    by (\ref{eq_proof_blow_up_K-_2}) again, if $T_2 \geq T_1$ is sufficiently large. Hence, one obtains
    \begin{equation}
        M^\prime(t)^2 \leq \frac{8}{9} M(t) M^{\prime \prime}(t), \quad t \geq T,
    \end{equation}
    implying that the function $M^{-\frac{1}{8}}$ is concave on $(T, + \infty)$, a contradiction with (\ref{eq_proof_blow_up_K-_2}).
\end{proof}

Finally, \eqref{eq_prop_duffing_K_energy_-infty} follows from the following lemma.

\begin{lem}\label{lem_blow-up_energy}
    Let $(u_0, u_1) \in \mathbb{R}^2$, $\gamma \geq 0$, and write $T^+ \in (0, +\infty]$ for the maximal time of existence of the solution $u$ of \eqref{duffing}. Assume that $T^+ < + \infty$. 
    Then
    \begin{equation}\label{eq_lem_blow-up_energy_1}
        \left\vert u(t) \right\vert \xrightarrow{t \rightarrow T^+} +\infty \quad \text{ and } \quad \left\vert \dot{u}(t) \right\vert \xrightarrow{t \rightarrow T^+} +\infty.
    \end{equation}
    If in addition $\gamma > 0$, then
    \begin{equation}\label{eq_lem_blow-up_energy_2}
         E_u(t) \xrightarrow{t \rightarrow T^+} -\infty.
    \end{equation}
\end{lem}

\begin{proof}
    First, we prove that, up to changing $u$ into $-u$, one has
    \begin{equation}\label{eq_proof_lem_blow-up_energy_1}
        u(t) \xrightarrow{t \rightarrow T^+} +\infty.
    \end{equation}
    Since $T^+ < + \infty$, the function $t \mapsto \left(u(t), \dot{u}(t) \right)$ is unbounded, and because $E_u(t) \leq E_u(0)$, this implies that $\vert u(t) \vert$ becomes arbitrarily large. Up to changing $u$ into $-u$, we can assume that there exists an increasing sequence $(t_n)$ such that $t_n$ tends to $T^+$ as $n$ tends to infinity, and 
    \begin{equation}
         u(t_n) \xrightarrow{n \rightarrow +\infty} +\infty.
    \end{equation}
    Up to a subsequence, we can assume that $\left( u(t_n) \right)$ is increasing. We use the following elementary lemma, whose proof is given below. 

    \begin{lem}\label{lem_elementary_f_prime_positive}
        Let $a < b$, and let $f \in \mathscr{C}^1([a, b], \mathbb{R})$ be such that $f(a) < f(b)$. Then, there exists $c \in [a, b]$ such that $f(c)>f(a)$ and $f^\prime(c) > 0$.
    \end{lem}

    By Lemma \ref{lem_elementary_f_prime_positive}, up to changing $t_n$ for some $s_n \in [t_n, t_{n+1}]$, we can assume that $\dot{u}(t_n) > 0$. Let $N$ be such that $u(t_N) > 1$. We claim that $\dot{u}(t) > 0$ for all $t \in [t_N, T^+)$. Assume by contradiction that there exists $T \in [t_N, T^+)$ such that $\dot{u}(T) = 0$ and $\dot{u}(t) > 0$ for all $t \in [t_N, T)$. For $t \in [t_N, T]$, one has $u(t) \geq u(t_N) > 1$, implying 
    \begin{equation}
        \ddot{u}(t) + \gamma \dot{u}(t) \geq u(t_N)^3 - u(t_N) > 0. 
    \end{equation}
    Multiplying by $e^{\gamma t}$ and integrating, one finds
    \begin{equation}
        \dot{u}(t) e^{\gamma t} \geq \dot{u}(t_N) e^{\gamma t_N} > 0, \quad t \in [t_N, T].
    \end{equation}
    For $t = T$, this gives a contradiction. Hence, $u$ is increasing on $[t_N, T^+)$, and thus \eqref{eq_proof_lem_blow-up_energy_1} holds true.
    
    Now, we complete the proof of \eqref{eq_lem_blow-up_energy_1}. Let $T_0 \in [0, T^+)$ be such that $u(t) > 1$ for all $t \in [T_0, T^+)$. There exists $T_1 \in [T_0, T^+)$, which can be chosen arbitrarily close to $T^+$, such that $\dot{u}(T_1) > 0$. As above, one has
    \begin{equation}
        \dot{u}(t) e^{\gamma t} \geq \dot{u}(T_1) e^{\gamma T_1} > 0, \quad t \in [T_1, T^+).
    \end{equation}
    Hence, one has
    \begin{equation}
        \dot{u}(s) = \sqrt{2 E_u(t) - u(s)^2 + \frac{u(s)^4}{2}} \leq \sqrt{2 E_u(0) - u(s)^2 + \frac{u(s)^4}{2}}, \quad s \in [T_1, T^+).
    \end{equation}
    Integrating and changing variables, this gives
    \begin{equation}
        \int_{u(t)}^{+ \infty} \frac{\dd x}{\sqrt{2 E_u(0) - x^2 + \frac{x^4}{2}}} \leq T^+ - t, \quad t \in [T_1, T^+). 
    \end{equation}
    By \eqref{eq_proof_lem_blow-up_energy_1}, up to increasing $T_1$, we can assume that if $t \in [T_1, T^+)$, then
    \begin{equation}
        \frac{1}{\sqrt{2 E_u(0) - x^2 + \frac{x^4}{2}}} \geq \frac{1}{x^2}, \quad x \geq u(t),
    \end{equation}
    which yields
    \begin{equation}
        u(t) \geq \frac{1}{T^+ - t}, \quad t \in [T_1, T^+).
    \end{equation}
    To complete the proof, we distinguish two cases. First, assume that the energy is bounded below, that is, that $E_u(t) \geq E_\infty$ for all $t \in [0, T^+)$, for some $E_\infty \in \mathbb{R}$. Then, there exists $T_2 \in [T_1, T^+)$ such that
    \begin{equation}\label{eq_proof_lem_blow-up_energy_2}
        \dot{u}(s) \geq \sqrt{2 E_\infty - u(s)^2 + \frac{u(s)^4}{2}} \geq \frac{c}{\left( T^+ - t \right)^2}, \quad t \in (T_2, T^+),
    \end{equation}
    for some $c>0$. This gives \eqref{eq_lem_blow-up_energy_1}. The fact that the energy is bounded implies that $\gamma = 0$. Indeed, if $\gamma > 0$, then using \eqref{eq_proof_lem_blow-up_energy_2} in \eqref{eq_energy_equality}, one finds \eqref{eq_lem_blow-up_energy_2}, a contradiction. Hence, $\gamma = 0$.
    
    Second, assume that the energy is not bounded below, implying in particular that $\gamma > 0$, and that \eqref{eq_lem_blow-up_energy_2} holds true. It remains to prove the second half of \eqref{eq_lem_blow-up_energy_1}. There exists $T_2 \in [T_1, T^+)$ such that
    \begin{equation}
        \ddot{u}(t) + \gamma \dot{u}(t) = u(t)^3 - u(t) \geq \frac{1}{2 \left( T^+ - t \right)^3}, \quad t \in (T_2, T^+). 
    \end{equation}
    Multiplying by $e^{\gamma t}$ and integrating, one obtains
    \begin{align}
        \dot{u}(t) e^{\gamma T^+} \geq \dot{u}(t) e^{\gamma t} & \geq \dot{u}(T_2) e^{\gamma T_2} + \int_{T_2}^t \frac{e^{\gamma \tau}}{2 \left( T^+ - \tau \right)^3} \dd \tau \\
        & \geq \dot{u}(T_2) e^{\gamma T_2} + \frac{e^{\gamma T_2}}{2} \int_{T_2}^t \frac{1}{\left( T^+ - \tau \right)^3} \dd \tau, \quad t \in (T_2, T^+). 
    \end{align}
    This proves that 
    \begin{equation}
        \dot{u}(s) \geq \frac{c}{\left( T^+ - t \right)^2}, \quad t \in (T_3, T^+),
    \end{equation}
    for some $c>0$ and $T_3 \in [T_2, T^+)$. This completes the proof.
\end{proof}

For completeness, we provide a short proof of Lemma \ref{lem_elementary_f_prime_positive}.

\begin{proof}[Proof of Lemma \ref{lem_elementary_f_prime_positive}.]
Let $x_0 = \sup \left\{ x \in [a,b],  f(x) \leq f(a) \right\}$. One has $x_0 < b$ and $f(x_0) = f(a)$. There exists $c \in (x_0, b)$ such that $f^\prime(c) = \frac{f(b) - f(x_0)}{b-x_0} > 0$. As $c > x_0$, one has $f(c) > f(a)$. 
\end{proof}

\paragraph{Solutions with initial data in $\mathcal{K}^+$.} Let $(u_0, u_1) \in \mathcal{K}^+$ and $\gamma \geq 0$. Arguing as above (see \eqref{eq:proof:potential_well}), one finds 
\begin{equation}\label{eq_proof_cv_K+_1}
    u(t)^2 < 1 - \delta, \quad t \geq 0,
\end{equation}
for some $\delta > 0$, implying also
\begin{equation}
    \frac{\dot{u}(t)^2}{2} \leq \frac{\dot{u}(t)^2}{2} + J(u(t)) = E_u(t) < \frac{1}{4}.
\end{equation}
Hence, the solution $u(t)$ is defined for all $t \geq 0$.

Now, we assume that $\gamma > 0$, and we prove that the solution tends to zero at an exponential rate. Note that (\ref{eq_proof_cv_K+_1}) implies
\begin{equation}\label{eq_proof_cv_K+_2}
    \frac{u(t)^2 + \dot{u}(t)^2}{4} \leq \frac{\dot{u}(t)^2}{2} + \frac{u(t)^2}{2} - \frac{u(t)^2 \times 1^2}{4} \leq E_u(t) \leq \frac{u(t)^2 + \dot{u}(t)^2}{2}, \quad t \geq 0.
\end{equation}
Let $\epsilon > 0$, and set $N = E_u + \epsilon u \dot{u}$. If $\epsilon$ is sufficiently small, then
\begin{equation}\label{eq_proof_cv_K+_3}
    \frac{1}{2} E_u(t) \leq N(t) \leq \frac{3}{2} E_u(t), \quad t \geq 0.
\end{equation}
One has
\begin{equation}
    N^\prime = (\epsilon - \gamma) \dot{u}^2 + \epsilon (u^4 - u^2) - \epsilon \gamma u \dot{u}. 
\end{equation}
Note that (\ref{eq_proof_cv_K+_1}) implies $u^4 - u^2 \leq - \delta u^2$. One has
\begin{equation}
    \epsilon \gamma \left\vert u \dot{u} \right\vert \leq C_0 \epsilon \sqrt{\epsilon} u^2 + \sqrt{\epsilon} \dot{u}^2,
\end{equation}
for some $C_0 > 0$. Hence, one obtains
\begin{equation}
    N^\prime \leq (\epsilon + \sqrt{\epsilon} - \gamma) \dot{u}^2 + \epsilon \left(C_0 \sqrt{\epsilon} - \delta \right) u^2, 
\end{equation}
implying that for $\epsilon$ sufficiently small, 
\begin{equation}
    N^\prime \leq - c_1 \left( \dot{u}^2 + u^2 \right), 
\end{equation}
for some $c_1 > 0$. By (\ref{eq_proof_cv_K+_2}) and (\ref{eq_proof_cv_K+_3}), this gives
\begin{equation}
    N^\prime \leq - c_2 N, 
\end{equation}
for some $c_2 > 0$, implying that $N$ tends to zero at an exponential rate. By (\ref{eq_proof_cv_K+_2}) and (\ref{eq_proof_cv_K+_3}) again, this proves that $(u, \dot{u})$ tends to zero at an exponential rate.

\subsection{Solutions of the Duffing equation of high energy}

Here, we prove Theorem \ref{thm_duffing_N}. Let $(u_0, u_1)$ be such that $E(u_0,u_1) \geq \frac{1}{4}$ with $u_1 > 0$, and let $\gamma \geq 0$.

\paragraph{Step 1: the case $E(u_0,u_1) = \frac{1}{4}$.} 
    Note that $u_1 > 0$ and $E(u_0,u_1) = \frac{1}{4}$ imply $u_0 \neq \pm 1$.
    If $\gamma > 0$, then as $u_1 > 0$ and $E^\prime_u(t) = - \gamma \dot{u}(t)^2$, one has $E_u(t) < \frac{1}{4}$ for all $t > 0$. In other words, the solution $u$ either enters $\mathcal{K}^+$ (if $u_0 \in (-1,1)$) and thus tends to zero at an exponential rate, or enters $\mathcal{K}^-$ (if $\vert u_0 \vert > 1$) and thus blows up.
    
    If $\gamma = 0$ and $E(u_0,u_1) = \frac{1}{4}$, then $u$ can be computed explicitly. Note that $\dot{u}(t) > 0$ for all $t \geq 0$, as $\dot{u}(t)=0$ and $E_u(t) = \frac{1}{4}$ would imply $\left( u(t), \dot{u}(t) \right) = (\pm 1, 0)$ for all $t \geq 0$. First, if $u_0 > 1$, then one has $u(t) \geq u_0 > 1$ for all $t \geq 0$, implying that $u$ blows up in finite time by Lemma \ref{lem_u_>_1_blowup}. Second, if $u_0 \in (-1, 1)$, then, using $E_u(t) = \frac{1}{4}$, one obtains
    \begin{equation}
        \dot{u}(t) = \frac{1 - u(t)^2}{\sqrt{2}}, \quad t \geq 0,
    \end{equation}
    and thus
    \begin{equation}
        u(t) = 1 - \frac{2(1 - u_0) e^{-t \sqrt{2}}}{1 + u_0 + (1 - u_0) e^{-t \sqrt{2}}}, \quad t \geq 0.
    \end{equation}
    In particular, $u$ tends to $1$ at an exponential rate. Third, arguing similarly, one proves that $u$ tends to $-1$ at an exponential rate if $u_0 < -1$.

\paragraph{Step 2: the case $E(u_0,u_1) > \frac{1}{4}$ and $\gamma = 0$.} 
    In that case, we prove that $u$ blows up in finite time. Assume that $u$ is defined on $[0, T]$, for some $T > 0$. One has
    \begin{equation}
        \frac{\dot{u}(t)^2}{2} = E_u(0) - J(u(t)) \geq E_u(0) - \frac{1}{4} > 0.
    \end{equation}
    In particular, if $T$ is sufficiently large, then there exist $\delta > 0$ and $T^\prime \geq 0$ such that $u(t) \geq 1 + \delta$ for all $t \geq T^\prime$. Hence, by Lemma \ref{lem_u_>_1_blowup}, $u$ blows up in finite time. 

\paragraph{Step 3: solutions with $E(u_0,u_1) > \frac{1}{4}$ and $\gamma > 0$, which do not exit $\mathcal{N}$.}
    Assume $\gamma > 0$ and $E_u(0) = E(u_0, u_1) > \frac{1}{4}$, and let $T^+ \in (0, +\infty]$ denote the maximal time of existence of $u$. Since $E_u$ is nonincreasing, we can distinguish two cases: first, there exists $T \in (0, T^+)$ such that $E_u(T) = \frac{1}{4}$, which brings us back to Step 1, or second, one has
    \begin{equation}
        E_u(t) \xrightarrow{t \rightarrow T^+} E_\infty \geq \frac{1}{4}.
    \end{equation}
    Assume that this second case occurs. By Lemma \ref{lem_blow-up_energy}, we know that $u$ is global. We prove that $E_\infty = \frac{1}{4}$ and that $u$ converges to $\pm 1$ at an exponential rate.
    
    Assume by contradiction that $E_\infty > \frac{1}{4}$. One has     
    \begin{equation}
        \frac{\dot{u}(t)^2}{2} \geq E_\infty - J(u(t)) \geq E_\infty - \frac{1}{4} > 0,
    \end{equation}
    implying that $\dot{u}$ does not vanish, and thus satisfies $\dot{u}(t) \geq \sqrt{2 E_\infty - \frac{1}{2}} > 0$. In particular, using Lemma \ref{lem_u_>_1_blowup} as above, one proves that $u$ blows up in finite time, a contradiction. Hence, $E_\infty = \frac{1}{4}$.
    
    Next, we prove that $u$ converges to $\pm 1$. If $\dot{u}(t) = 0$, then $E_u(t) \leq \frac{1}{4}$, with equality if and only if $u(t) = \pm 1$. Hence, one has $E_u(t) > \frac{1}{4}$ for all $t \in [0, T^+)$. In particular, $\dot{u}$ does not vanish, and therefore $\dot{u}(t) > 0$ for all $t \geq 0$. If there exists $t \geq 0$ such that $u(t) > 1$, then Lemma \ref{lem_blow-up_geq_1_plus_delta} implies that $u$ blows up in finite time, a contradiction. Hence, $u$ is bounded, implying in particular that $u$ converges to some $u_\infty \in \mathbb{R}$. Since $E_u$ converges, this implies that $\dot{u}$ also converges to some limit. Since $u$ is bounded, one obtains
    \begin{equation}
        \dot{u}(t) \xrightarrow{t \rightarrow +\infty} 0.
    \end{equation}
    Letting $t$ tend to infinity in the definition of $E_u$, one finds $J(u_\infty) = \frac{1}{4}$, and thus $u_\infty = \pm 1$.

    Finally, we prove that $(u, \dot{u})$ converges to $(\pm 1, 0)$ at an exponential rate. We assume that $u$ converges to $1$, the other case being similar. In particular, one has $u(t) \leq 1$ for all $t \geq 0$. Set $\Tilde{E}(t) = E_u(t) - \frac{1}{4} \geq 0$, and $M(t) = \left( u(t) - 1 \right)^2$. One has 
    \begin{equation}
        \Tilde{E}(t) = \frac{\dot{u}(t)^2}{2} - \frac{\left( u(t) - 1 \right)^2 \left( u(t) + 1 \right)^2}{4}.
    \end{equation}
    Since $u$ converges to $1$, there exists $T > 0$ such that 
    \begin{equation}
        \frac{1}{2} \leq \frac{\left( u(t) + 1 \right)^2}{4} \leq \frac{3}{2}, \quad t \geq T,
    \end{equation}
    implying
    \begin{equation}\label{eq_proof_solution_not_exit_N_1}
        \frac{\dot{u}(t)^2}{2} - \frac{3 \left( u(t) - 1 \right)^2}{2} \leq \Tilde{E}(t) \leq \frac{\dot{u}(t)^2}{2} - \frac{\left( u(t) - 1 \right)^2}{2} \leq \frac{\dot{u}(t)^2}{2}, \quad t \geq T.
    \end{equation}
    In particular, this yields
    \begin{equation}
        \Tilde{E}^\prime(t) = - \gamma \dot{u}(t)^2 \leq - 2 \gamma \Tilde{E}(t), \quad t \geq T,
    \end{equation}
    and therefore
    \begin{equation}\label{eq_proof_solution_not_exit_N_2}
        0 \leq \Tilde{E}(t) \leq \Tilde{E}(T) e^{- 2 \gamma(t - T)}, \quad t\geq T.
    \end{equation}
    Using \eqref{eq_proof_solution_not_exit_N_1} and $\Tilde{E}\geq 0$, one finds
    \begin{equation}
        M^\prime(t) = 2 \dot{u}(t) \left( u(t) - 1 \right) = - 2 \left\vert \dot{u}(t) \right\vert \left\vert u(t) - 1 \right\vert \leq - 2 \left( u(t) - 1 \right)^2, \quad t \geq T,
    \end{equation}
    and therefore
    \begin{equation}
        0 \leq M(t) \leq M(T) e^{- 2 (t - T)}, \quad t \geq T.
    \end{equation}
    Hence, $u$ converges to $1$ at an exponential rate. By \eqref{eq_proof_solution_not_exit_N_1} and \eqref{eq_proof_solution_not_exit_N_2}, this implies that $\dot{u}$ converges to $0$ at an exponential rate.
    
\begin{figure}[ht]
        \centering
        \begin{tikzpicture}[scale=2.5,>=latex]
            \draw[->] (-1.4,0) -- (1.4,0) node[right] {$u(t)$};
            \draw[->] (0,-1.4) -- (0,1.4) node[above] {$\dot{u}(t)$};
            
            \draw[dashed,domain=-1.6:1.6,smooth,thick,variable=\x] plot ({\x},{(\x*\x - 1)/sqrt(2)});
            \draw[dashed,domain=-1.6:1.6,smooth,thick,variable=\x] plot ({\x},{-(\x*\x - 1)/sqrt(2)});
            
            \draw[thick,gray]
            (-0.5,0.78) 
              .. controls +(0.3,0.17) and +(-0.3,0.14) .. (0.42,0.92)
              .. controls +(0.3,-0.14) and +(-0.15,0.0) .. (1.0,0.68)
              .. controls +(0.15,0.0) and +(-0.036,-0.06) .. (1.3,0.95);

            \draw[thick,black]
            (-0.5,0.78) 
              .. controls +(0.3,0.11) and +(-0.3,0.17) .. (0.42,0.76)
              .. controls +(0.3,-0.17) and +(-0.15,0.0) .. (1.0,0.38)
              .. controls +(0.15,0.0) and +(-0.036,-0.06) .. (1.3,0.7);

            \fill (1,0) circle(0.03);
            \fill (-1,0) circle(0.03);
            \node at (1,-0.2) {$1$};
            \node at (-1,-0.2) {$-1$};

            \fill (-0.5,0.78) circle(0.03);
            \node at (-0.75,0.71) {$(u_0, u_1)$};

            \draw[densely dotted] (0.42,0) -- (0.42,0.92);
            \draw[densely dotted] (0.42,0.76) -- (1.8,0.76);
            \draw[densely dotted] (0.42,0.92) -- (1.8,0.92);

            \fill (0.42,0.76) circle(0.03);
            
            \fill[gray] (0.42,0.92) circle(0.03);
            
            \draw[<->] (1.8,0.76) -- (1.8,0.92);

            \node at (0.42,-0.1) {$x$};
            \node at (2.01,0.84) {$\phi(x)$};

            \node[gray] at (-0.1,1.04) {$u^1$};
            \node at (-0.1,0.79) {$u^2$};

        \end{tikzpicture}
        \caption{The function $\phi$ of Lemma \ref{lem_comparison_solutions}, with the solution $u^1$ corresponding to damping $\gamma_1$ (in gray) and the solution $u^2$ corresponding to damping $\gamma_2 > \gamma_1$ (in black).}
        \label{fig_step_4_comparaison}
    \end{figure}
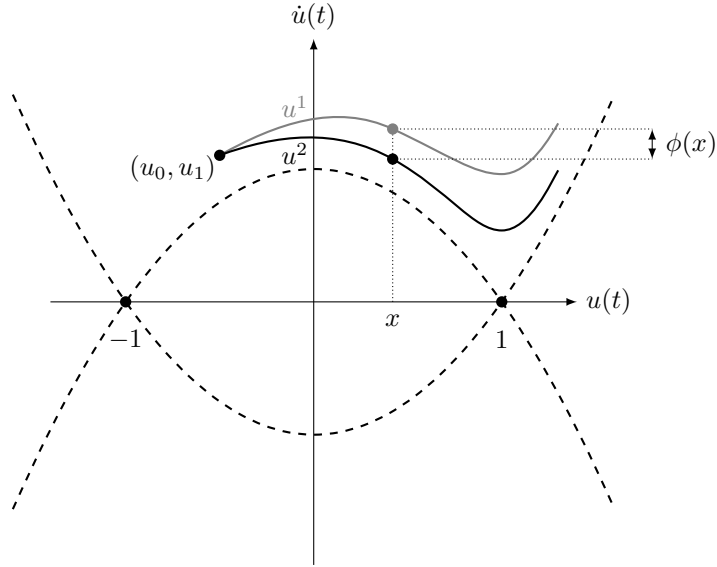

\paragraph{Step 4: comparison between solutions with different dampings.}
    We prove the following technical result, which will be used repeatedly in what follows. An illustrative picture is shown in Figure~\ref{fig_step_4_comparaison}.
    
    \begin{lem}\label{lem_comparison_solutions}
        Let $\gamma_2 > \gamma_1 \geq 0$, and $T_2, T_1 \in (0, +\infty]$. Let $u^1$ and $u^2$ be solutions of \eqref{duffing} with damping $\gamma_1$ and $\gamma_2$, respectively, with the same initial data $(u_0, u_1)$. For $j \in \{1,2\}$, assume that $u^j$ is defined on $[0, T_j)$ and that $(u^j)^\prime(t) > 0$ for all $t \in [0, T_j)$. Set $I = [u_0, u^1(T_1)] \cap [u_0, u^2(T_2)]$, and
        \begin{equation}
            \phi: x \in I \longmapsto \dot{u}^1\left( (u^1)^{-1}(x) \right) - \dot{u}^2\left( (u^2)^{-1}(x) \right).
        \end{equation}
        Then for all $x \in I \setminus \{ u_0 \}$, one has $\phi(x) > 0$. In addition, for all $x \in I$, if $x \in (- \infty, -1] \cup [0, 1]$, then $\phi^\prime(x) > 0$.
    \end{lem}

    \begin{proof}
        Note that the functions $(u^1)^{-1}$ and $(u^2)^{-1}$ are well-defined and continuous on $I$, and that $\phi(u_0)=0$. One has
        \begin{equation}\label{eq_proof_lem_comparison_solutions_1}
            \phi^\prime(x) = \frac{(x-x^3) \phi(x)}{\dot{u}^1\left( (u^1)^{-1}(x) \right)  \dot{u}^2\left( (u^2)^{-1}(x) \right) } + \gamma_2 - \gamma_1, \quad x \in I.
        \end{equation}
        As $\phi^\prime(u_0) > 0$, one has $\phi(u_0 + \epsilon) > 0$ if $\epsilon>0$ is sufficiently small. Assume by contradiction that $\phi(x_0)=0$ for some $x_0 \in I$, with $x_0 > 0$. By the intermediate value theorem, up to reducing $x_0$, we can assume that $\phi(x) > 0$ for all $x \in (0, x_0)$. Since $\phi^\prime(x_0) = \gamma_2 - \gamma_1 > 0$, one has $\phi(x_0 - \epsilon) < 0$ if $\epsilon > 0$ is sufficiently small, a contradiction. 

        If $x \in (- \infty, -1] \cup [0, 1]$, then $x - x^3 \geq 0$, so that $\phi^\prime(x) > 0$ directly follows from \eqref{eq_proof_lem_comparison_solutions_1}.
    \end{proof}

    We are now ready to study the behaviour of solutions with initial data in $\mathcal{N}_1$, $\mathcal{N}_2$, and $\mathcal{N}_3$.

\paragraph{Step 5: the case $(u_0, u_1) \in \mathcal{N}_1$.} 
    Assume $\gamma > 0$, $u_0 \geq 1$, $u_1 > 0$ and $E(u_0, u_1) > \frac{1}{4}$. Write $T^+ > 0$ for the maximal time of existence of $u$. We prove that for all $t \in [0, T^+)$, $\dot{u}(t) > 0$. Assume by contradiction that there exists $T \in [0, T^+)$ such that $\dot{u}(T) = 0$ and $\dot{u}(t) > 0$ for all $t \in [0, T)$. One has $u(T) > u(0) \geq 1$, yielding
    \begin{equation}
        \ddot{u}(T) = u(T)^3 - u(T) > 0.
    \end{equation}
    Hence, there exists $\epsilon > 0$ such that $\ddot{u}(t) > 0$ for all $t \in [T-\epsilon, T]$. In particular, one has $\dot{u}(T) > \dot{u}(T-\epsilon) > 0$, a contradiction. It implies that $u(t) \geq u\left( \frac{T^+}{2} \right) > 1$ for all $t \in \left[ \frac{T^+}{2}, T^+\right)$, and therefore $u$ blows up by Lemma \ref{lem_blow-up_geq_1_plus_delta}.

\paragraph{Step 6: the case $(u_0, u_1) \in \mathcal{N}_2$.} 
    Consider $1 > u_0 \geq -1$ and $u_1 > 0$ such that $E(u_0, u_1) > \frac{1}{4}$, and for $\gamma \geq 0$, write $T^+_\gamma > 0$ for the maximal time of existence of the solution $u$ with initial data $(u_0, u_1)$ and with damping $\gamma$. Set
    \begin{equation}
        A = \left\{ \gamma \geq 0, T^+_\gamma < + \infty \right\}.
    \end{equation}
    By Step 2, we know that $0 \in A$. By Step 3, if $\gamma \geq 0$ does not belong to $A$, then either $E_u(t) > \frac{1}{4}$ for all $t \geq 0$ and 
    \begin{equation}\label{eq_proof_N_2_def_B}
        \left( u(t), \partial_t u(t) \right) \xrightarrow{t \rightarrow + \infty} (1,0),
    \end{equation}
    or there exists $T > 0$ such that $E_u(T) = \frac{1}{4}$ and $u(T) \in (-1, 1)$, and hence  
    \begin{equation}\label{eq_proof_N_2_def_C}
        \left( u(t), \partial_t u(t) \right) \xrightarrow{t \rightarrow + \infty} (0,0).
    \end{equation}
    Hence, we define
    \begin{equation}
        B = \left\{ \gamma \geq 0, T^+_\gamma = + \infty \text{ and \eqref{eq_proof_N_2_def_B} occurs} \right\} \quad \text{and} \quad C = \left\{ \gamma \geq 0, T^+_\gamma = + \infty \text{ and \eqref{eq_proof_N_2_def_C} occurs} \right\},
    \end{equation}
    so that $[0, +\infty) = A \sqcup B \sqcup C$. We prove that there exists $\gamma_0 > 0$ such that $A = [0, \gamma_0)$, $B = \{\gamma_0\}$, and $C = (\gamma_0, +\infty)$.
    
    \underline{First, we prove that $A$ is bounded.}
    To improve the readability of this argument, we isolate the key point in a lemma. 
    
    \begin{lem}\label{lem_technic_bound_on_gamma}
        Let $\gamma > 0$ and $T_2 > T_1 \geq 0$. Assume that $u$ is a solution of \eqref{duffing} on $[T_1, T_2]$ such that $u(t) \in [0, 1]$ for $t \in [T_1, T_2]$ and $\dot{u}(T_1) > 0$ and $u(T_2) = 1$. Then
        \begin{equation}\label{eq_lem_technic_bound_on_gamma}
            1 \leq u(T_1) + \frac{\dot{u}(T_1)}{\gamma}.
        \end{equation}
    \end{lem}

    \begin{proof}
        Since $u(t) \in [0, 1]$ for $t \in [T_1, T_2]$, one has $\ddot{u} + \gamma \dot{u} \leq 0$ on $[T_1, T_2]$, implying
        \begin{equation}
            \dot{u}(t) \leq \dot{u}(T_1) e^{\gamma ( T_1 - t )}, \quad t \in [T_1, T_2].
        \end{equation}
        Integrating, and using $\dot{u}(T_1) > 0$, one finds
        \begin{equation}
            u(t) \leq u(T_1) + \dot{u}(T_1) \frac{1 - e^{\gamma ( T_1 - t )}}{\gamma} \leq u(T_1) + \frac{\dot{u}(T_1)}{\gamma}, \quad t \in [T_1, T_2].
        \end{equation}
        For $t = T_2$, this gives \eqref{eq_lem_technic_bound_on_gamma}.
    \end{proof}
    
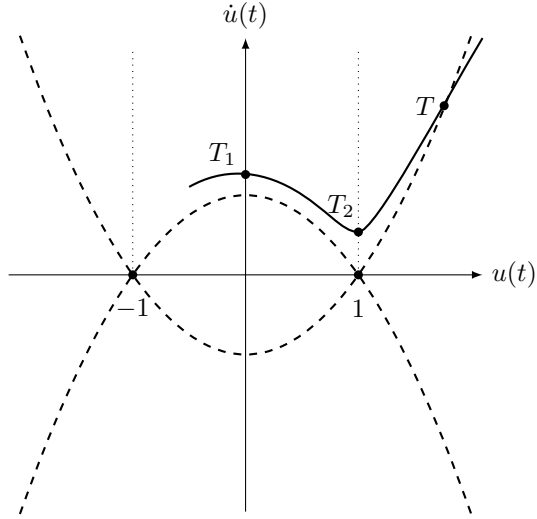
\begin{figure}[ht]
        \centering
        \begin{tikzpicture}[scale=1.5,>=latex]
            \draw[->] (-2.1,0) -- (2.1,0) node[right] {$u(t)$};
            \draw[->] (0,-2.1) -- (0,2.1) node[above] {$\dot{u}(t)$};
            
            \draw[dotted] (1,0) -- (1,2);
            \draw[dotted] (-1,0) -- (-1,2);
            
            \draw[dashed,domain=-2:2,smooth,thick,variable=\x] plot ({\x},{(\x*\x - 1)/sqrt(2)});
            \draw[dashed,domain=-2:2,smooth,thick,variable=\x] plot ({\x},{-(\x*\x - 1)/sqrt(2)});
            
            \draw[thick]
            (-0.5,0.78) 
              .. controls +(0.3,0.17) and +(-0.3,0.17) .. (0.42,0.76)
              .. controls +(0.3,-0.17) and +(-0.15,0.0) .. (1.0,0.38)
              .. controls +(0.15,0.0) and +(-0.36,-0.6) .. (2.1,2.1);
            
            \fill (1,0) circle(0.04);
            \fill (-1,0) circle(0.04);
            \node at (1,-0.3) {$1$};
            \node at (-1,-0.3) {$-1$};
            
            \fill (0,0.89) circle(0.04);
            \node[above left] at (0,0.89) {$T_1$};
            \fill (1,0.38) circle(0.04);
            \node[above left] at (1.05,0.43) {$T_2$};
            \fill (1.76,1.5) circle(0.04);
            \node[left] at (1.76,1.5) {$T$};
        \end{tikzpicture}
        \caption{The times $T$, $T_1$ and $T_2$ in the proof of the boundedness of $A$ (in the case $u_0 < 0$).}
        \label{fig_step_6_A_bounded}
    \end{figure}

    Now, let $\gamma \in A$, with $\gamma > 0$. Since $T^+_\gamma < + \infty$, Step 3 implies the existence of $T \in (0, T^+)$ such that $E_u(T) = \frac{1}{4}$ and $E_u(t) > \frac{1}{4}$ for $t \in [0, T)$. This implies that $\dot{u}(t) > 0 $ for $t \in [0, T)$. Since $T^+_\gamma < + \infty$, Step 1 yields $u(T) > 1$. In particular, $(u(T), \partial_t u(T)) \in \mathcal{N}_1$. By the intermediate value theorem, there exists $T_2 \in [0, T]$ such that $u(T_2) = 1$ and $u(t) \leq 1$ for all $t \in [0, T_2]$. Finally, we define $T_1$ depending on the sign of $u_0$: if $u_0 \geq 0$, we simply set $T_1 = 0$; if $u_0 < 0$, we apply the intermediate value theorem again to define $T_1 \in [0, T_2]$ such that $u(T_1) = 0$ and $u(t) \in [0, 1]$ for all $t \in [T_1, T_2]$. The situation is illustrated in Figure \ref{fig_step_6_A_bounded}. Now, Lemma \ref{lem_technic_bound_on_gamma} gives
    \begin{equation}
        1 \leq u(T_1) + \frac{\dot{u}(T_1)}{\gamma}.
    \end{equation}
    If $T_1 = 0$, this gives $\gamma \leq \frac{u_1}{1 - u_0}$. If $u(T_1) = 0$, then this gives 
    \begin{equation}
        \gamma \leq \dot{u}(T_1) = \sqrt{2 E_u(T_1)} \leq \sqrt{2 E(u_0, u_1)} .
    \end{equation}
    Hence, we have proved that
    \begin{equation}
        A \subset \left[0, \max \left( \frac{u_1}{1 - u_0}, \sqrt{2 E(u_0, u_1)} \right) \right].
    \end{equation}

    \underline{Second, we prove that $A$ is an interval.} Let $\gamma_2 \in A$, with $\gamma_2 > 0$, and let $\gamma_1 \in (0, \gamma_2)$. We prove that $\gamma_1 \in A$. Let $u^1$ and $u^2$ be the solutions of \eqref{duffing} with damping $\gamma_1$ and $\gamma_2$ respectively. By Step 3, we know that either there exists $T_1 \in (0, T^+_{\gamma_1})$ such that $E_{u^1}(T_1) = \frac{1}{4}$, or 
    \begin{equation}\label{eq_proof_N_2_1}
        T^+_{\gamma_1} = +\infty \quad \text{ and } \quad \left( u^1(t), \partial_t u^1(t) \right) \xrightarrow{t \rightarrow + \infty} (1,0).
    \end{equation}
    
    We prove that the latter case cannot occur. Assume by contradiction that \eqref{eq_proof_N_2_1} holds true. Then $\dot{u}^1(t) > 0$ for all $t \geq 0$. We know that there exists $T_2 \in (0, T^+_{\gamma_2})$ such that $u^2(T_2)=1$ and $\dot{u}^2(t) > 0$ for all $t \in [0, T_2]$. Hence, Lemma \ref{lem_comparison_solutions} gives 
    \begin{equation}
        \dot{u}^1(t) > \dot{u}^2\left( (u^2)^{-1}\left( u^1(t) \right) \right), \quad t \geq 0.
    \end{equation}
    Letting $t$ tend to $+ \infty$ gives $\dot{u}^2(T_2) \leq 0$, a contradiction. This proves that \eqref{eq_proof_N_2_1} cannot occur.

    Hence, we have proved that there exists $T_1 \in (0, T^+_{\gamma_1})$ such that $E_{u^1}(T_1) = \frac{1}{4}$. Assume by contradiction that $u^1(T_1) \leq 1$. Set $t_2 = (u^2)^{-1}(u^1(T_1))$, so that $u^2(t_2) = u^1(T_1) \leq 1$. Then, Lemma \ref{lem_comparison_solutions} gives $\dot{u}^1(T_1) > \dot{u}^2(t_2)$, implying 
    \begin{equation}
        E_{u^2}(t_2) = J(u^1(T_1)) + \frac{\dot{u}^2(t_2)^2}{2} \leq E_{u^1}(T_1) = \frac{1}{4}. 
    \end{equation} 
    Hence, by Step 1, $T^+_{\gamma_2} = +\infty$, a contradiction. Hence, $u^1(T_1) > 1$ and thus $\gamma_1 \in A$.

    Summarising, we have proved that there exists $\gamma_0 \geq 0$ such that 
    \begin{equation}
        \overline{A} = [0, \gamma_0].
    \end{equation}

    \underline{Third, we prove that $\gamma_0 \notin A$ and that $C$ is open.} Let $\gamma \in A$, and write $u$ for the solution of \eqref{duffing} with damping $\gamma$. There exists $T \in (0, T^+_{\gamma})$ such that $u(T) > 1$ and $\dot{u}(T) > 0$. By standard ODE theory, for all $\epsilon > 0$, there exists $\delta > 0$ such that for all $\gamma^\prime \in [\gamma - \delta, \gamma + \delta]$, the solution $v$ of \eqref{duffing} with damping $\gamma^\prime$ is defined on $[0, T]$ and satisfies
    \begin{equation}
        \left\vert u(T) - v(T) \right\vert + \left\vert \dot{u}(T) - \dot{v}(T) \right\vert \leq \epsilon.  
    \end{equation}
    Hence, if $\gamma^\prime \geq 0$ is sufficiently close to $\gamma$, then $\gamma^\prime \in A$. In particular, $\gamma_0 > 0$ and 
    \begin{equation}
        A = [0, \gamma_0).
    \end{equation}
    For all $\gamma > 0$, a solution with initial data close to $(0,0)$ is in $\mathcal{K}^+$ and thus converges to zero at an exponential rate. Hence, a similar argument shows that $C$ is open. 

    \underline{Fourth, we prove that if $\gamma_1 \in C$ and $\gamma_2 > \gamma_1$ then $\gamma_2 \in C$.} Let $u^1$ and $u^2$ be the solutions of \eqref{duffing} with damping $\gamma_1$ and $\gamma_2$ respectively. Since $\gamma_2 \notin A$, one has $\gamma_2 \in B \sqcup C$. Assume by contradiction that $\gamma_2 \in B$. Let $T_1 > 0$ be such that $E_{u^1}(T_1) = \frac{1}{4}$. One has $\dot{u}^1 > 0$ on $[0, T_1]$, $\dot{u}^2 > 0$ on $[0, +\infty)$, and $u^1(T_1) < 1$. Set $T_2 = (u^2)^{-1}(u^1(T_1)) > 0$. As above, one has $\dot{u}^2(T_2) < \dot{u}^1(T_1)$, and 
    \begin{equation}
        E_{u^2}(T_2) = J(u^1(T_1)) + \frac{\dot{u}^2(T_2)^2}{2} < E_{u^1}(T_1) = \frac{1}{4},
    \end{equation}
    a contradiction. Hence, $\gamma_2 \in C$.

    \underline{Fifth, we prove that if $\gamma_1 \in B$ and $\gamma_2 > \gamma_1$ then $\gamma_2 \in C$.} Let $u^1$ and $u^2$ be the solutions of \eqref{duffing} with damping $\gamma_1$ and $\gamma_2$ respectively. Assume by contradiction that $\gamma_2 \in B$. One has $\dot{u}^1 > 0$ and $\dot{u}^2 > 0$ on $[0, +\infty)$. Set
    \begin{equation}
        \phi: x \in [u_0, 1) \longmapsto \dot{u}^1\left( (u^1)^{-1}(x) \right) - \dot{u}^2\left( (u^2)^{-1}(x) \right)
    \end{equation}
    By Lemma \ref{lem_comparison_solutions}, one has $\phi(x) \geq 0$ for all $x \in [u_0, 1)$, and $\phi^\prime(x) > 0$ for all $x \in [u_0, 1) \cap (0, 1)$. By assumption, $\phi$ converges to $0$ as $x$ tends to $1$, a contradiction. Summarising, we have proved that $C = (\gamma_0, +\infty)$ and $B = \{\gamma_0\}$.

\paragraph{Step 7: the case $(u_0, u_1) \in \mathcal{N}_3$.}
    Many arguments are analogous to those in Step 6, so we omit some repetitive details and provide full explanations only for the more delicate points. Consider $u_0 < -1$ and $u_1 > 0$ such that $E(u_0, u_1) > \frac{1}{4}$, and for $\gamma \geq 0$, write $T^+_\gamma > 0$ for the maximal time of existence of the solution $u$ with initial data $(u_0, u_1)$ and with damping $\gamma$. Set
    \begin{equation}
        \begin{array}{c}
            A = \left\{\gamma \geq 0, \text{there exists } T \in (0, T^+_\gamma) \text{ such that } u(T) \geq -1 \right\},\\
            B = \left\{\gamma \geq 0, u(t) < -1 \text{ for all } t \in [0, T^+_\gamma) \text{ and } T^+_\gamma = + \infty\right\},\\
            C = \left\{\gamma \geq 0, u(t) < -1 \text{ for all } t \in [0, T^+_\gamma) \text{ and } T^+_\gamma < + \infty\right\}.
        \end{array}
    \end{equation}
    One has $[0, +\infty) = A \sqcup B \sqcup C$, and $0 \in A$. For all $\gamma \geq 0$, note that: first, if $\gamma \in A$ then there exists $t \in (0, T^+_\gamma)$ such that $u(t) = -1$, implying $E_u(t) \geq \frac{1}{4}$ and therefore $\left( u(t), \partial_t u(t) \right) \in \mathcal{N}_2$; second, if $\gamma \in B$ then $E_u(t) > \frac{1}{4}$ for all $t \geq 0$ (otherwise $u$ would enter $\mathcal{K}^-$ and therefore blow up in finite time) implying that $\left( u(t), \partial_t u(t) \right)$ converges to $-1$ at an exponential rate; third, if $\gamma \in C$ then $u$ enters $\mathcal{K}^-$.

    As above (see Lemma \ref{lem_technic_bound_on_gamma}), one proves that if $T>0$ is such that $u(t) \leq -1$ for $t \in [0, T]$ and $u(T) = -1$, then
    \begin{equation}
        -1 \leq u_0 + \frac{u_1}{\gamma},
    \end{equation}
    implying that $A$ is bounded. As before, one proves that $A$ is an interval. If $\gamma \in A$, then one proves that there exists $T \in (0, T^+_\gamma)$ such that $u(T) > -1$ and $\dot{u}(T) > 0$, and, as above, this implies that $A$ is open. Hence, there exists $\gamma_1 > 0$ such that $A = [0, \gamma_1)$. Arguing as before, one proves that $C = (\gamma_1, + \infty)$ and $B = \{\gamma_1\}$.

    Next, we further decompose $A$ into $A = [0, \gamma_1) = A^\prime \sqcup B^\prime \sqcup C^\prime$, with
    \begin{equation}
        \begin{array}{c}
            A^\prime = \left\{ \gamma \in A, T^+_\gamma < + \infty \right\},\\
            B^\prime = \left\{ \gamma \in A, T^+_\gamma = + \infty \text{ and }
            \left( u(t), \partial_t u(t) \right) \xrightarrow{t \rightarrow + \infty} (1,0) \right\},\\
            C^\prime = \left\{ \gamma \in A, T^+_\gamma = + \infty \text{ and }
            \left( u(t), \partial_t u(t) \right) \xrightarrow{t \rightarrow + \infty} (0,0) \right\}.
        \end{array}
    \end{equation}

    As above, one has the following properties: $0 \in A^\prime$, $C^\prime$ is an open interval, there exists $\gamma_0 \in (0, \gamma_1]$ such that $A^\prime = [0, \gamma_0)$, $B^\prime$ contains at most one element, $(\gamma, \gamma^\prime) \in A^\prime \times B^\prime$ implies $\gamma < \gamma^\prime$, and $(\gamma, \gamma^\prime) \in B^\prime \times C^\prime$ implies $\gamma < \gamma^\prime$. To complete the proof, we need to show that $\gamma_0 < \gamma_1$, or equivalently, that $C^\prime$ (and hence $B^\prime$) is nonempty. The key argument is given in the following lemma, whose proof is given below. Intuitively, the idea is the following: for a fixed value of $\gamma > 0$, if $(u_0, u_1)$ is close to $(-1, 0)$, with $u_1 > 0$, then the solution enters $\mathcal{K}^+$.
    
    \begin{lem}\label{lem_u_1_close_to_-1}
        There exists an absolute constant $c_0 > 0$ satisfying the following property. For all $\gamma > 0$ and all $(u_0, u_1) \in \mathbb{R}^2$ such that $u_0 \in \left[ -1, -\frac{1}{2} \right]$ and $u_1 \in \left(0, 2 \gamma \right)$, if there exists $T > 0$ such that $u(T) = 0$ and
        \begin{equation}
            E_u(t) > \frac{1}{4} \quad t \in [0, T],
        \end{equation}
        then one must have $u_1 > c_0 \min(\gamma, 1)$.
    \end{lem}

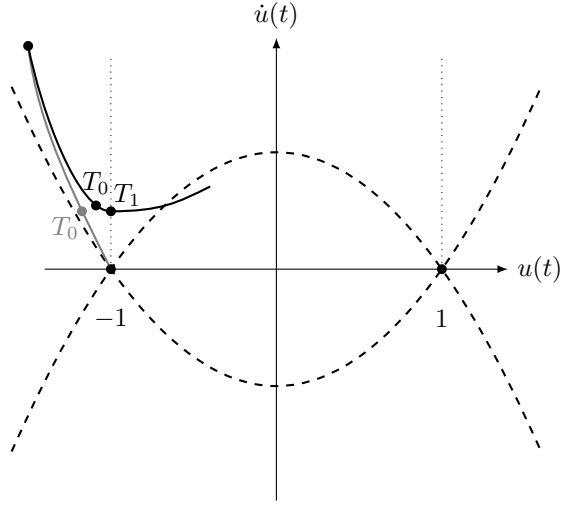
\begin{figure}[ht]
        \centering
        \begin{tikzpicture}[scale=2.2,>=latex]
            \draw[->] (-1.4,0) -- (1.4,0) node[right] {$u(t)$};
            \draw[->] (0,-1.4) -- (0,1.4) node[above] {$\dot{u}(t)$};
            
            \draw[dotted] (1,0) -- (1,1.3);
            \draw[dotted] (-1,0) -- (-1,1.3);
            
            \draw[dashed,domain=-1.6:1.6,smooth,thick,variable=\x] plot ({\x},{(\x*\x - 1)/sqrt(2)});
            \draw[dashed,domain=-1.6:1.6,smooth,thick,variable=\x] plot ({\x},{-(\x*\x - 1)/sqrt(2)});
            
            \draw[thick, gray]
            (-1.5,1.35) .. controls +(0.06,-0.5) and +(-0.2,0.4) .. (-1,0);

            \draw[thick]
            (-1.5,1.35) 
              .. controls +(0.08,-0.4) and +(-0.2,0) .. (-1,0.35)
              .. controls +(0.3,0) and +(-0.2,-0.1) .. (-0.4,0.5);
            
            \fill (1,0) circle(0.03);
            \fill (-1,0) circle(0.03);
            \node at (1,-0.3) {$1$};
            \node at (-1,-0.3) {$-1$};

            \fill (-1.5,1.35) circle(0.03);
            \fill (-1,0.35) circle(0.03);
            \node at (-0.9,0.45) {$T_1$};
            \fill (-1.09,0.385) circle(0.03);
            \node at (-1.09,0.52) {$T_0$};
            \fill[gray] (-1.175,0.35) circle(0.03);
            \node[gray] at (-1.275,0.25) {$T_0$};
        \end{tikzpicture}
        \caption{The times $T_0$ and $T_1$ in the proof of $C^\prime \neq \emptyset$ (in the case $u_\gamma(T_0) < -1$). The curves in black and gray represent $u_{\gamma}$ and $u_{\gamma_1}$ respectively.}
        \label{fig_step_7_C_prime_nonempty}
    \end{figure}

    We use Lemma \ref{lem_u_1_close_to_-1} to prove that $C^\prime$ is nonempty. In this paragraph, we write $u_\gamma$ for the solution of \eqref{duffing} with damping $\gamma \geq 0$ (and with initial data $(u_0, u_1)$, defined in the beginning of Step 7). Assume that $\gamma_1 \leq 1$, so that $\min(\gamma_1, 1) = \gamma_1$, the other case being similar. Let $\epsilon > 0$ be such that
    \begin{equation}\label{eq_proof_step_7_1}
        \epsilon < \frac{1}{2}, \quad \epsilon < \gamma_1, \quad \text{and} \quad \epsilon < \frac{c_0 \gamma_1}{2}.
    \end{equation}
    Since $\left( u_{\gamma_1}, \dot{u}_{\gamma_1} \right)$ converges to $(-1,0)$, there exists $T_0 > 0$ such that 
    \begin{equation}
        \left\vert u_{\gamma_1}(T_0) + 1 \right\vert + \left\vert \dot{u}_{\gamma_1}(T_0) \right\vert \leq \frac{\epsilon}{2}.
    \end{equation}
    By standard ODE theory, there exists $\delta = \delta(\epsilon) > 0$ such that for all $\gamma \geq 0$, if $\vert \gamma_1 - \gamma \vert \leq \delta$ then $T^+_\gamma > T_0$ and
    \begin{equation}
        \left\vert u_{\gamma_1}(T_0) - u_\gamma(T_0) \right\vert + \left\vert \dot{u}_{\gamma_1}(T_0) - \dot{u}_\gamma(T_0) \right\vert \leq \frac{\epsilon}{2},
    \end{equation}
    implying
    \begin{equation}
        \left\vert u_{\gamma}(T_0) + 1 \right\vert + \left\vert \dot{u}_{\gamma}(T_0) \right\vert \leq \epsilon.
    \end{equation}
    Let $\gamma \in \left(\gamma_1 - \delta, \gamma_1 \right)$, with $\gamma \geq 0$. Write $u = u_\gamma$. Since $\gamma \in A$, there exists $T_1 > 0$ such that $u(T_1) = -1$. Hence, we know that either $(u, \dot{u})$ converges (to $(0, 0)$ or to $(1, 0)$), or one has $E_u(t) > \frac{1}{4}$ for all $t \geq 0$ such that $u(t) \leq 1$. We prove that the latter case cannot occur, implying that $\gamma \notin A^\prime$, and hence $C^\prime$ and $B^\prime$ are nonempty. Up to reducing $\delta$, \eqref{eq_proof_step_7_1} implies $\epsilon < 2 \gamma$, and $\epsilon < c_0 \gamma$. If $u(T_0) \geq -1$, then we can apply Lemma \ref{lem_u_1_close_to_-1} with the initial data $\left( \Tilde{u}_0, \Tilde{u}_1 \right) = \left( u(T_1), \dot{u}(T_1) \right)$: it implies that one cannot have $E_u(t) > \frac{1}{4}$ for all $t$ such that $u(t) \in [-1, 0]$, yielding $\gamma \notin A^\prime$. If $u(T_0) < -1$, then $T_0 < T_1$, and one has $\ddot{u} < 0$ on $(T_0, T_1)$, implying $0 < \dot{u}(T_1) < \dot{u}(T_0) \leq \epsilon$. The situation is illustrated in Figure \ref{fig_step_7_C_prime_nonempty}. Hence, we can apply Lemma \ref{lem_u_1_close_to_-1} with the initial data $\left( \Tilde{u}_0, \Tilde{u}_1 \right) = \left( u(T_0), \dot{u}(T_0) \right)$, to get the same conclusion. 

    \begin{proof}[Proof of Lemma \ref{lem_u_1_close_to_-1}.]
        Let $\gamma > 0$, $T > 0$ and $(u_0, u_1) \in \left[-1,-\frac{1}{2}\right] \times \left(0, 2 \gamma \right)$, be such that $u(T) = 0$ and $E_u(t) > \frac{1}{4}$ for all $t \in [0, T]$. In particular, $\dot{u}(t) > 0$ for all $t \in [0, T]$. The energy equality yields
        \begin{equation}
            \frac{u_1^2}{2} > \frac{1}{4} + \frac{u_1^2}{2} - E_u(T) \geq E_u(0) - E_u(T) = \gamma \int_0^T \dot{u}(s)^2 \dd s,
        \end{equation}
        and it gives 
        \begin{equation}
            \frac{1}{2} \leq -u_0 = \int_0^T \dot{u}(s) \dd s \leq \sqrt{T} \left( \int_0^T \dot{u}(s)^2 \dd s \right)^{\frac{1}{2}} \leq \frac{u_1 \sqrt{T}}{\sqrt{2 \gamma}}.
        \end{equation}
        By the mean value theorem, there exists $t \in [0, T]$ such that $\dot{u}(t) = \frac{u(T) - u(0)}{T} \leq \frac{1}{T}$, implying
        \begin{equation}\label{eq_proof_lem_u_1_close_to_-1_1}
            u_1 \geq \sqrt{\frac{\gamma}{2} \min_{[0,T]} \dot{u}}.
        \end{equation}
        We claim that
        \begin{equation}\label{eq_proof_lem_u_1_close_to_-1_claim}
            \min_{[0, T]} \dot{u} \geq \min\left( \frac{u_1}{2}, \frac{u_1}{2 \sqrt{2} \gamma} \right).
        \end{equation}
        If $\gamma \leq \frac{1}{\sqrt{2}}$, then \eqref{eq_proof_lem_u_1_close_to_-1_1} and \eqref{eq_proof_lem_u_1_close_to_-1_claim} give $u_1 \geq \frac{\gamma}{4}$, and if $\gamma > \frac{1}{\sqrt{2}}$, then \eqref{eq_proof_lem_u_1_close_to_-1_1} and \eqref{eq_proof_lem_u_1_close_to_-1_claim} give $u_1 \geq \frac{1}{4 \sqrt{2}}$, completing the proof.

        It remains to prove \eqref{eq_proof_lem_u_1_close_to_-1_claim}. Let $t_0 \in [0, T]$ be such that $\dot{u}(t_0) = \min_{[0, T]} \dot{u}$. The idea is the following: if $u(t_0)$ is close to $-1$, then $\dot{u}(t_0)$ is close to $u_1$, whereas if $u(t_0)$ is close to $0$, the condition $E_u(t_0) \geq \frac{1}{4}$ prevents $\dot{u}(t_0)$ from being too small. We must balance these two effects to establish \eqref{eq_proof_lem_u_1_close_to_-1_claim}. On the one hand, one has
        \begin{equation}
            \dot{u}(t_0) = u_1 + \int_0^{t_0} \left( u(s)^3 - u(s) - \gamma \dot{u}(s) \right) \dd s \geq u_1 - \gamma u(t_0) - \gamma.
        \end{equation}
        On the other hand, one has
        \begin{equation}
            \frac{\dot{u}(t_0)^2}{2} = E_u(t_0) - \frac{u(t_0)^2}{2} + \frac{u(t_0)^4}{4} \geq \frac{1}{4} - \frac{u(t_0)^2}{2} + \frac{u(t_0)^4}{4}.
        \end{equation}
        Hence, \eqref{eq_proof_lem_u_1_close_to_-1_claim} follows from the following elementary lemma.

        \begin{lem}\label{lem_technical_X_Y_min_dot_u}
            Let $\gamma > 0$ and $u_1 \in (0, 2 \gamma)$. Let $X \in [-1, 0]$ and $Y > 0$. Assume
            \begin{equation}\label{eq_lem_technical_X_Y_min_dot_u}
                \frac{Y^2}{2} \geq \frac{1}{4} - \frac{X^2}{2} + \frac{X^4}{4} \quad \text{ and } \quad Y \geq u_1 - \gamma X - \gamma.
            \end{equation}
            Then $Y \geq \min\left( \frac{u_1}{2}, \frac{u_1}{2 \sqrt{2} \gamma} \right)$.
        \end{lem}

        \begin{proof}[Proof of Lemma \ref{lem_technical_X_Y_min_dot_u}.]
            Set $X_0 = \frac{u_1}{2\gamma} - 1 \in (-1, 0)$. If $X \leq X_0$, then
            \begin{equation}
                Y \geq u_1 - \gamma X - \gamma \geq \frac{u_1}{2}.
            \end{equation}
            If $X > X_0$, then 
            \begin{equation}
                Y \geq \frac{1 - X^2}{\sqrt{2}} \geq \frac{1 + X}{\sqrt{2}} \geq \frac{1 + X_0}{\sqrt{2}} = \frac{u_1}{2 \sqrt{2} \gamma}.
            \end{equation}            
        \end{proof}

        The proof of Lemma \ref{lem_u_1_close_to_-1} is complete. 
    \end{proof}

    Summarising, we have proved that there exists $0 < \gamma_0 < \gamma_1$ such that $[0, +\infty) = A^\prime \sqcup B^\prime \sqcup C^\prime \sqcup B \sqcup C$, with $A^\prime = [0, \gamma_0)$, $B^\prime = \{\gamma_0\}$, $C^\prime = (\gamma_0, \gamma_1)$, $B = \{\gamma_1\}$, and $C = (\gamma_1, + \infty)$. This completes the proof of Theorem \ref{thm_duffing_N}.
 
\subsection{Additional result about the Duffing equation}
 
\begin{figure}[ht]
    \centering
    \includegraphics[width=0.3\textwidth]{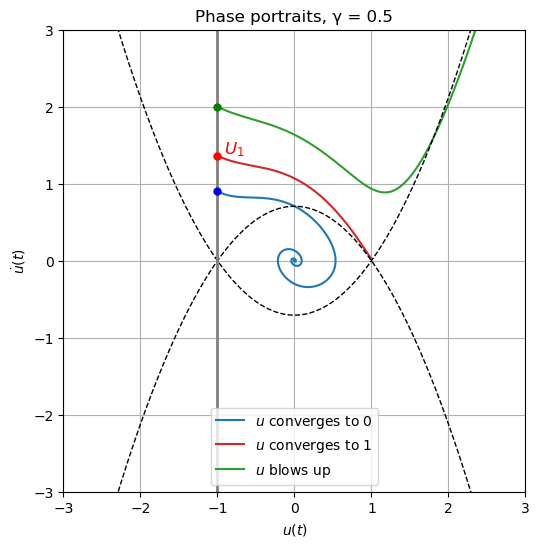}
    \caption{Illustration of Lemma~\ref{lem:gamma_fix_u_1_varies}: for $u_0 = -1$ and a fixed damping $\gamma > 0$, varying the initial velocity $u_1$ leads to different behaviours of the solution (convergence to $(0,0)$, convergence to $(1,0)$, or finite-time blow-up).}
    \label{fig_illustration_lem_gamma_constant}
\end{figure}

To construct solutions \eqref{KG} with interesting behaviour, we will need the following lemma concerning solutions of \eqref{duffing}. The key difference between this result and Theorem~\ref{thm_duffing_N} is that, in Theorem~\ref{thm_duffing_N}, we fix $(u_0, u_1)$ and study the behaviour of the solution as $\gamma$ varies, whereas here we fix $\gamma > 0$ and aim to prove the existence of initial data leading to different types of behaviour. An illustration of Lemma \ref{lem:gamma_fix_u_1_varies} is shown in Figure \ref{fig_illustration_lem_gamma_constant}.

\begin{lem}\label{lem:gamma_fix_u_1_varies}
    Consider $\gamma \in (0, +\infty)$. There exists $U_1 > 0$ such that for all $u_1 > 0$, if $u$ denotes the solution of \eqref{duffing} with initial data $(-1, u_1)$ and damping $\gamma$, then:
    \begin{enumerate}[label=(\roman*)]
        \item if $u_1 \in (0, U_1)$, then 
            \begin{equation}\label{eq:lem:gamma_fix_u_1_varies_1}
                \left( u(t), \dot{u}(t) \right) \xrightarrow{t \rightarrow + \infty} (0,0),
            \end{equation}
            at an exponential rate,
         \item if $u_1 = U_1$, then 
            \begin{equation}
                \left( u(t), \dot{u}(t) \right) \xrightarrow{t \rightarrow + \infty} (1,0),
            \end{equation}
            at an exponential rate,
        \item if $u_1 > U_1$, then $u$ blows up in finite time. 
    \end{enumerate}
\end{lem}

\begin{proof}
Note that if $u_1 = 0$, then $u(t) = -1$ for all $t \geq 0$. Write $T^+(u_1) \in (0, +\infty]$ for the maximal time of existence of the solution of \eqref{duffing} with initial data $(-1, u_1)$ and damping $\gamma$, and set
\begin{equation}
    \begin{split}
        & A = \left\{ u_1 > 0, T^+(u_1) = + \infty \text{ and } \left(u(t), \dot{u}(t) \right) \xrightarrow{t \rightarrow + \infty} (0,0) \right\} , \\
        & B = \left\{ u_1 > 0, T^+(u_1) = + \infty \text{ and } \left(u(t), \dot{u}(t) \right) \xrightarrow{t \rightarrow + \infty} (1,0) \right\} , \\
        & C = \left\{ u_1 > 0, T^+(u_1) < + \infty \right\} .
    \end{split}
\end{equation}
By Theorem \ref{thm_duffing_N}, we know that $(0, +\infty) = A \sqcup B \sqcup C$, that solutions with $u_1 \in A \sqcup B$ satisfy $u(t) \in [-1, 1)$ for all $t \geq 0$, while solutions with $u_1 \in C$ satisfy $u(t) \geq 1$ for some $t \in \left[ 0, T^+(u_1) \right)$. 

\paragraph{Step 1.}
Let $u_1 \in A \sqcup B$ and let $u$ be the solution of \eqref{duffing} with initial data $(-1, u_1)$ and damping $\gamma$. We prove that there exists $C_0 =  C_0(\gamma) > 0$ such that $u_1 \leq C_0$, implying that $A \sqcup B$ is bounded, and hence that $C$ is nonempty. Set $c = \frac{2}{3 \sqrt{3}}$. Since $u(t) \in [-1, 1)$ for all $t \geq 0$, one has 
\begin{equation}
    \ddot{u}(t) + \gamma \dot{u}(t) = u(t)^3 - u(t) \geq - c. 
\end{equation}
Let $v$ be the solution of $\ddot{v}(t) + \gamma \dot{v}(t) = - c$, with $\left( v(0), \dot{v}(0) \right) = \left( u(0), \dot{u}(0) \right) = (-1, u_1)$, and set $w = u - v$. Then, for all $t \geq 0$, one has
\begin{equation}
    \ddot{w}(t) + \gamma \dot{w}(t) \geq 0,
\end{equation}
implying that $\dot{w}(t) e^{\gamma t} \geq 0$, and hence $w(t) \geq 0$, that is, $u(t) \geq v(t)$. In addition, for all $t \geq 0$, a straightforward integration yields
\begin{equation}
    \dot{v}(t) = u_1 e^{-\gamma t} + \frac{c}{\gamma} \left( e^{-\gamma t} - 1 \right),
\end{equation}
and
\begin{equation}
    v(t) = - 1 + \frac{u_1}{\gamma} \left( 1 - e^{-\gamma t} \right) + \frac{c}{\gamma} \left( \frac{1 - e^{-\gamma t}}{\gamma} - t \right).
\end{equation}
Set $t_0 = \frac{\ln(2)}{\gamma}$. One has 
\begin{equation}
    u(t_0) \geq v(t_0) = -1 + \frac{u_1}{2 \gamma} + \frac{c}{\gamma} \left( \frac{1}{2\gamma} - \frac{\ln(2)}{\gamma} \right).
\end{equation}
Hence, there exists $C_0 =  C_0(\gamma) > 0$ such that if $u_1 > C_0$, then $u(t_0) \geq 1$, implying $u_1 \in C$, a contradiction. Hence, $u_1 \leq C_0$.

\paragraph{Step 2.} We prove that if $u_1 > 0$ is small, then $u_1 \in A$. Let $c_0 > 0$ be given by Lemma \ref{lem_u_1_close_to_-1}. Set $c_0^\prime = \min \left(c_0, c_0 \gamma, 2 \gamma \right)$. Let $u_1 \in (0, c_0^\prime)$, and let $u$ be the solution of \eqref{duffing} with initial data $(-1, u_1)$ and damping $\gamma$. If $u_1 \in B$, then one has $E_u(t) > \frac{1}{4}$ for all $t \geq 0$, and $u(t) \rightarrow 1$. If $u_1 \in C$, then there exists $T > 0$ such that $u(t) = 1$ and $E_u(t) > \frac{1}{4}$ for all $t \in [0, T]$. Hence, if $u_1 \in B \sqcup C$, then by the intermediate value theorem, there exists $T \in \left[ 0, T^+(u_1) \right)$ such that $u(T) = 0$ and $E_u(t) > \frac{1}{4}$ for all $t \in [0, T]$: Lemma \ref{lem_u_1_close_to_-1} thus gives $u_1 > \min(c_0, c_0 \gamma)$, a contradiction. Hence, $u_1 \in A$, and this proves that $(0, c_0^\prime) \subset A$.

\paragraph{Step 3.} As in the proof of Theorem \ref{thm_duffing_N}, we can compare solutions. More precisely, let $u_1^1 > u_1^2 > 0$, and write $u^1$ for the solution with initial data $(-1, u_1^1)$ and damping $\gamma$, and $u^2$ for the solution with initial data $(-1, u_1^2)$ and damping $\gamma$. For $j \in \{1,2\}$, let $T_j > 0$ be such that $\dot{u}^j(t) > 0$ for all $t \in [0, T_j]$. Arguing as in the proof of Lemma \ref{lem_comparison_solutions}, one proves that
\begin{equation}
    \dot{u}^1\left( (u^1)^{-1}(x) \right) > \dot{u}^2\left( (u^2)^{-1}(x) \right), 
\end{equation}
for all $x \in I := [-1, u^1(T_1)] \cap [-1, u^2(T_2)]$. As in the proof of Theorem \ref{thm_duffing_N}, it implies that if $u_1^1 \in A$, then $u_1^2 \in A$; whereas if $u_1^2 \in B \sqcup C$, then $u_1^1 \in C$. Hence, $A$ and $C$ are intervals, and $B$ contains at most one element. As in the proof of Theorem \ref{thm_duffing_N}, standard ODE theory shows that $A$ and $C$ are open. Therefore, $B$ is non-empty, and there exists $U_1$ such that $B = \{U_1\}$, $A = (0, U_1)$ and $C = (U_1, +\infty)$.
\end{proof}
 
\section{The cubic Klein-Gordon equation on a compact manifold}\label{sec:KG}

Recall that $\Omega$ is defined at the beginning of the introduction. We will assume that $\partial \Omega = \emptyset$ only when proving Theorem \ref{thm:construction_solution_particuliere} and Proposition \ref{prop:Q_nonconstant_lambda_1_smaller_than_2}.

\subsection{Well-posedness and continuity estimates for the cubic Klein-Gordon equation}

The following result is standard and we omit its proof. Note that in the case $T < + \infty$, the blow-up of $u(t)$ in $H_0^1(\Omega)$ follows from the energy estimate and the embedding $H^1(\Omega) \hookrightarrow L^4(\Omega)$ (see \cite[Theorem 11]{perrin2024damped}).

\begin{thm}\label{thm_existence_damped_waves}
    Consider $\gamma \in L^\infty(\Omega, \mathbb{R}_+)$. For any (real-valued) initial data $\left( u_0, u_1 \right) \in H_0^1(\Omega) \times L^2(\Omega)$, there exist a maximal time of existence $T \in (0, + \infty]$ and a unique solution $u$ of \eqref{KG} in $\mathscr{C}^0([0, T), H_0^1(\Omega)) \cap \mathscr{C}^1([0, T), L^2(\Omega))$. If $T < + \infty$, then 
    \begin{equation}\label{eq_thm_existence_wave_blow_up}
        \left\Vert u(t) \right\Vert_{H_0^1} \xrightarrow{t \rightarrow T^-} + \infty.
    \end{equation}
\end{thm}

We will need a source-to-solution continuity result for solutions of \eqref{KG}. Although it relies on standard techniques, we include a proof for completeness.

\begin{lem}\label{lem:source_to_sol}
    Let $T > 0$, $\gamma \in L^\infty(\Omega, \mathbb{R}_+)$, and $\left( u_0, u_1 \right) \in H_0^1(\Omega) \times L^2(\Omega)$. Assume that the solution $u$ of \eqref{KG} is defined on $[0, T]$. Then, there exist $C > 0$ and a neighbourhood $\mathcal{O}$ of $\left( u_0, u_1 \right) \in H_0^1(\Omega) \times L^2(\Omega)$ such that for all $\left( v_0, v_1 \right) \in \mathcal{O}$, the solution $v$ of \eqref{KG} with initial data $\left( v_0, v_1 \right)$ and damping $\gamma$ is defined on $[0,T]$ and satisfies
    \begin{equation}
        \sup_{t \in [0, T]} \left\Vert \left( u, \partial_t u \right)(t) -  \left( v, \partial_t v \right)(t) \right\Vert_{H_0^1(\Omega) \times L^2(\Omega)} \leq C \left\Vert \left( u_0, u_1 \right) -  \left( v_0, v_1 \right) \right\Vert_{H_0^1(\Omega) \times L^2(\Omega)} .
    \end{equation}
\end{lem}

\begin{proof}
    We split the proof into two steps. For $T > 0$, write $X_T = \mathscr{C}^0([0, T), H_0^1(\Omega)) \cap \mathscr{C}^1([0, T), L^2(\Omega))$ and
    \begin{equation}
        \Vert u \Vert_{X_T} = \sup_{t \in [0, T]} \left\Vert \left( u, \partial_t u \right)(t) \right\Vert_{H_0^1(\Omega) \times L^2(\Omega)}.
    \end{equation}
    
    \paragraph{Step 1: Solutions of the linear equation with time-dependent potential.} Consider $T > 0$, $\gamma \in L^\infty(\Omega, \mathbb{R}_+)$, and $\left( u_0, u_1 \right) \in H_0^1(\Omega) \times L^2(\Omega)$, $V \in L^2((0, T), L^3(\Omega))$ and $F \in L^1((0, T), L^2(\Omega))$. We prove that the solution $u$ of
    \begin{equation}\label{eq:proof:lem:source_to_sol_1}
        \left \{
            \begin{array}{rcccl}
                \partial_t^2 u - \Delta u + \gamma \partial_t u + V u & = & F & \quad & \text{in } (0, T) \times \Omega, \\
                (u(0), \partial_t u(0)) & = & \left( u_0, u_1 \right) & \quad & \text{in } \Omega, \\
                u & = & 0 & \quad & \text{in } (0, T) \times \partial \Omega,
            \end{array}
        \right.
    \end{equation}
    satisfies $u \in X_T$ and
    \begin{equation}\label{eq:proof:lem:source_to_sol_2}
        \Vert u \Vert_{X_T} \leq C_1 \left( \left\Vert \left( u_0, u_1 \right) \right\Vert_{H_0^1(\Omega) \times L^2(\Omega)} + \left\Vert F \right\Vert_{L^1((0, T), L^2)} \right),
    \end{equation}
    for some $C_1 > 0$ which is independent of $\left( u_0, u_1 \right)$ and $F$. The case $V = 0$ is well known; we omit its proof and use it to derive the result for $V \neq 0$. We may assume that $T \leq \delta$, where $\delta > 0$ is a small constant which is independent of $\left( u_0, u_1 \right)$ and $F$; the result for larger $T$ then follows by iteration. 
    
    For $U \in X_T$, let $u = \Gamma(U)$ be the solution of 
    \begin{equation}
        \left \{
            \begin{array}{rcccl}
                \partial_t^2 u - \Delta u + \gamma \partial_t u & = & F - VU & \quad & \text{in } (0, T) \times \Omega, \\
                (u(0), \partial_t u(0)) & = & \left( u_0, u_1 \right) & \quad & \text{in } \Omega, \\
                u & = & 0 & \quad & \text{in } (0, T) \times \partial \Omega.
            \end{array}
        \right.
    \end{equation}
    The case $V = 0$ gives 
    \begin{equation}
        \left\Vert \Gamma(U) \right\Vert_{X_T} \leq C_1 \left( \left\Vert \left( u_0, u_1 \right) \right\Vert_{H_0^1(\Omega) \times L^2(\Omega)} + \left\Vert F \right\Vert_{L^1((0, T), L^2)} + \left\Vert VU \right\Vert_{L^1((0, T), L^2)} \right).
    \end{equation}
    Using Hölder's inequality and the Sobolev embedding $H^1(\Omega) \hookrightarrow L^6(\Omega)$, one finds
    \begin{equation}
         \left\Vert VU \right\Vert_{L^1((0, T), L^2)} \leq \left\Vert V \right\Vert_{L^1((0, T), L^3)} \left\Vert U \right\Vert_{L^\infty((0, T), H^1)} \leq C_2 \sqrt{T} \left\Vert V \right\Vert_{L^2((0, T), L^3)} \left\Vert U \right\Vert_{X_T},
    \end{equation}
    for some $C_2 > 0$ (which depends only on $\Omega$). Similarly, one has
    \begin{equation}
        \left\Vert \Gamma(U) - \Gamma\left(\widetilde{U}\right) \right\Vert_{X_T} \leq C_1 \left\Vert V \left( U - \widetilde{U}\right) \right\Vert_{L^1((0, T), L^2)} \leq C_1 C_2 \sqrt{T} \left\Vert V \right\Vert_{L^2((0, T), L^3)} \left\Vert U - \widetilde{U} \right\Vert_{X_T}.
    \end{equation}
    Assume that $C_1 C_2 \sqrt{T} \left\Vert V \right\Vert_{L^2((0, T), L^3)} \leq \frac{1}{2}$, and set $R = 2 C_1 \left( \left\Vert \left( u_0, u_1 \right) \right\Vert_{H_0^1(\Omega) \times L^2(\Omega)} + \left\Vert F \right\Vert_{L^1((0, T), L^2)} \right)$. Then for $U$, $\widetilde{U} \in X_T$ such that $\Vert U \Vert_{X_T} \leq R$, one has
    \begin{equation}
        \left\Vert \Gamma(U) \right\Vert_{X_T} \leq R \quad \text{ and } \left\Vert \Gamma(U) - \Gamma\left(\widetilde{U}\right) \right\Vert_{X_T} \leq \frac{1}{2} \left\Vert U - \widetilde{U} \right\Vert_{X_T}.
    \end{equation}
    The Picard fixed point theorem proves that the solution $u$ of \eqref{eq:proof:lem:source_to_sol_1} exists, and satisfies $\left\Vert u \right\Vert_{X_T} \leq R$. This completes the proof of \eqref{eq:proof:lem:source_to_sol_2} for small $T$.

    \paragraph{Step 2: Source-to-solution continuity.} Let $T > 0$, $\gamma \in L^\infty(\Omega, \mathbb{R}_+)$, and $u$ be as in the statement of Lemma \ref{lem:source_to_sol}. The idea is to apply the Picard fixed point theorem to $h = v - u$. More precisely, proving the existence of $v$ is equivalent to proving the existence of $h$ satisfying
    \begin{equation}
        \left \{
            \begin{array}{rcccl}
                \partial_t^2 h - \Delta h + \gamma \partial_t h + (1 - 3 u^2) h & = & 3 u h^2 + h^3 & \quad & \text{in } (0, T) \times \Omega, \\
                (h(0), \partial_t h(0)) & = & \left( v_0 - u_0, v_1 - u_1 \right) & \quad & \text{in } \Omega, \\
                h & = & 0 & \quad & \text{in } (0, T) \times \partial \Omega.
            \end{array}
        \right.
    \end{equation}
    We reuse the notation $\Gamma$ to denote a different mapping from that defined in Step 1. Set $\left( h_0, h_1 \right) = \left( v_0 - u_0, v_1 - u_1 \right)$, and for $H \in X_T$, let $h = \Gamma(H)$ be the solution of 
    \begin{equation}
        \left \{
            \begin{array}{rcccl}
                \partial_t^2 h - \Delta h + \gamma \partial_t h + (1 - 3 u^2) h & = & 3 u H^2 + H^3 & \quad & \text{in } (0, T) \times \Omega, \\
                (h(0), \partial_t h(0)) & = & \left( h_0, h_1 \right) & \quad & \text{in } \Omega, \\
                h & = & 0 & \quad & \text{in } (0, T) \times \partial \Omega.
            \end{array}
        \right.
    \end{equation}
    Note that the potential $V := 1 - 3 u^2$ satisfies $V \in L^2((0, T), L^3)$ by the Sobolev embedding $H^1(\Omega) \hookrightarrow L^6(\Omega)$. Let $C_1 > 0$ be the constant given by Step 1 (which may depend on $T$, $u$, and $\gamma$). One has     
    \begin{equation}
        \left\Vert \Gamma(H) \right\Vert_{X_T} \leq C_1 \left( \left\Vert \left( h_0, h_1 \right) \right\Vert_{H_0^1(\Omega) \times L^2(\Omega)} + \left\Vert 3 u H^2 + H^3 \right\Vert_{L^1((0, T), L^2)} \right).
    \end{equation}
    Using Hölder's inequality and the Sobolev embedding $H^1(\Omega) \hookrightarrow L^6(\Omega)$, one finds
    \begin{equation}
        \left\Vert 3 u H^2 + H^3 \right\Vert_{L^1((0, T), L^2)} \leq C_3 \left( \Vert H \Vert_{X_T}^2 + \Vert H \Vert_{X_T}^3 \right), 
    \end{equation}
    for some $C_3> 0$ which may depend on $u$ and $T$. Similarly, one proves that
    \begin{equation}
        \left\Vert \Gamma(H) - \Gamma\left( \widetilde{H}\right) \right\Vert_{X_T} \leq C_4 \left\Vert H - \widetilde{H} \right\Vert_{X_T} \left( \left\Vert H \right\Vert_{X_T} + \left\Vert H \right\Vert_{X_T}^2 + \left\Vert \widetilde{H} \right\Vert_{X_T} + \left\Vert \widetilde{H} \right\Vert_{X_T}^2 \right),
    \end{equation}
    for some $C_4> 0$ which may depend on $u$ and $T$. Set $R = 2 C_1 \left\Vert \left( h_0, h_1 \right) \right\Vert_{H_0^1(\Omega) \times L^2(\Omega)}$. If $R$ is sufficiently small, one has
    \begin{equation}
        C_1 C_3 \left( R^2 + R^3 \right) \leq \frac{R}{2} \quad \text{ and } \quad 2 C_4 \left( R + R^2 \right) \leq \frac{1}{2},
    \end{equation}
    so that the conclusion follows from the Picard fixed point theorem. Hence, the open set $\mathcal{O}$ in Lemma \ref{lem:source_to_sol} can be chosen as a small ball in $H_0^1(\Omega) \times L^2(\Omega)$, centered at $(u_0, u_1)$, with a radius depending on $u$ and $T$. 
\end{proof}

\subsection{Ground states of the cubic Klein-Gordon equation}

We recall the definition and some properties of ground states. Recall that $K$ is defined by \eqref{eq_def_K}.

\begin{lem}\label{lem_ground_state_properties}
    Write 
    \begin{equation}
        \mathcal{N}^\prime = \left\{ w \in H_0^1(\Omega), w \neq 0, -\Delta w + \beta w = w^3 \right\} \quad \text{ and } \quad \mathcal{N} = \left\{ w \in H_0^1(\Omega), w \neq 0, K(w) = 0 \right\},
    \end{equation}
    so that by definition, $d = \inf \left\{ J_\KG(w), w \in \mathcal{N}^\prime \right\}$. The set $\mathcal{N}$ is called the Nehari manifold. One has $d = \inf \left\{ J_\KG(w), w \in \mathcal{N} \right\}$ and $d > 0$. Any solution $Q \in H_0^1(\Omega)$ of $-\Delta Q + \beta Q = Q^3$ such that $J_\KG(Q) = d$ is called a ground state of \eqref{KG}. There exists (at least) one ground state, and if $Q$ is a ground state, then $Q \in H^2(\Omega)$.
\end{lem}

\begin{rem}
    If $Q$ is a ground state, then so is $-Q$. The uniqueness of positive ground states is, in general, false (see the bibliographic comments in the introduction).
\end{rem}

\begin{proof}
    To simplify the notation, we write $J = J_\KG$ throughout this proof.
    
    First, we prove that $d > 0$. Note that $\mathcal{N}^\prime \subset \mathcal{N}$. Let $w \in \mathcal{N}$. One has
    \begin{equation}
        \left\Vert w \right\Vert_{H_0^1}^2 = \left\Vert w \right\Vert_{L^4}^4 + K(w) = \left\Vert w \right\Vert_{L^4}^4 \leq C \left\Vert w \right\Vert_{H_0^1}^4,
    \end{equation}
    for some $C>0$, by the Sobolev embedding $H^1(\Omega) \hookrightarrow L^4(\Omega)$. Hence, one obtains
    \begin{equation}
        J(w) = \frac{1}{4} \left\Vert w \right\Vert_{H_0^1}^2 + \frac{K(w)}{4} = \frac{1}{4} \left\Vert w \right\Vert_{H_0^1}^2 \geq \frac{1}{4C} > 0,
    \end{equation}
    and therefore $d \geq \inf \left\{J(w), w \in \mathcal{N} \right\} > 0$.

    Second, we prove the existence of a ground state, and we prove that $d = \inf_{\mathcal{N}} J$. Consider a sequence $\left( w_n \right)_{n \in \mathbb{N}}$ of elements of $\mathcal{N}$ such that $J(w_n) \longrightarrow \inf_{\mathcal{N}} J$ as $n$ tends to infinity. Up to a subsequence, we may assume that there exists $w \in H_0^1(\Omega)$ such that $w_n$ converges to $w$ weakly in $H^1(\Omega)$ and strongly in $L^4(\Omega)$. One has
    \begin{align}
        \inf_{\mathcal{N}} J & = \lim_{n \rightarrow \infty} \left( \frac{1}{4} \Vert w_n \Vert_{H_0^1}^2 + \frac{K(w_n)}{4}\right) = \frac{1}{4} \lim_{n \rightarrow \infty} \Vert w_n \Vert_{H_0^1}^2 \\
        & = \lim_{n \rightarrow \infty} \left( \frac{1}{4} \Vert w_n \Vert_{H_0^1}^2 - \frac{K(w_n)}{4}\right) = \frac{1}{4} \lim_{n \rightarrow \infty} \Vert w_n \Vert_{L^4}^4 = \frac{1}{4} \Vert w \Vert_{L^4}^4,
    \end{align}
    implying in particular that $w \neq 0$. As $w_n$ converges to $w$ weakly, one has $\Vert w \Vert_{H_0^1} \leq \lim \Vert w_n \Vert_{H_0^1}$, implying $K(w) \leq 0$. For $\lambda > 0$, write $j(\lambda) = J(\lambda w)$, and let $\lambda^\ast > 0$ be the argument of the unique maximum of $j$. One has $K(\lambda^\ast w)=j^\prime(\lambda^\ast) =0$, implying
    \begin{equation}
        \inf_{\mathcal{N}} J \leq J(\lambda^\ast w) \leq \left( \frac{(\lambda^\ast)^2}{2} - \frac{(\lambda^\ast)^4}{4} \right) \lim_{n \rightarrow \infty} \Vert w_n \Vert_{H_0^1}^2 \leq \frac{1}{4} \lim_{n \rightarrow \infty} \Vert w_n \Vert_{H_0^1}^2 = \inf_{\mathcal{N}} J.
    \end{equation}
    implying that the inequality $\frac{(\lambda^\ast)^2}{2} - \frac{(\lambda^\ast)^4}{4} \leq \frac{1}{4}$ is an equality, and therefore $\lambda^\ast = 1$ and $K(w) = 0$. Hence, one obtains $\Vert w \Vert_{H_0^1} = \lim \Vert w_n \Vert_{H_0^1}$ and $J(w) = \inf_{\mathcal{N}} J$. We claim that $w \in \mathcal{N}^\prime$, implying $d = \inf_{\mathcal{N}} J$. 

    As $K(w) = 0$, one has
    \begin{equation}\label{eq_proof_lem_min_d_implies_min_J}
        \left\langle \Delta w + \beta w - w^3, w \right\rangle_{H^{-1}(\Omega) \times H_0^1(\Omega)} = 0.
    \end{equation}
    The differential of $J$ and $K$ are given by
    \begin{equation}
        \dd J(u)(h) = \int_\Omega \left( \nabla u \cdot \nabla h + \beta u h \right) \dd x - \int_\Omega u^3 h \dd x, \quad u, h \in H_0^1(\Omega),
    \end{equation}
    and
    \begin{equation}
        \dd K(u)(h) = 2 \int_\Omega \left( \nabla u \cdot \nabla h + \beta u h \right) \dd x - 4 \int_\Omega u^3 h \dd x \quad u, h \in H_0^1(\Omega).
    \end{equation}
    Hence, there exists a Lagrange multiplier $\Lambda \in \mathbb{R}$ such that the Euler equation
    \begin{equation}
        \left( - \Delta w + \beta w - w^3 \right) + \Lambda \left( - 2 \Delta w + 2 \beta w - 4 w^3 \right) = 0
    \end{equation}
    holds in $H^{-1}(\Omega)$. Rewriting that equation, one finds
    \begin{equation}
        \left(1 + 2 \Lambda \right) \left( - \Delta w + \beta w - w^3 \right) - 2 \Lambda w^3 = 0.
    \end{equation}
    By (\ref{eq_proof_lem_min_d_implies_min_J}), this implies 
    \begin{equation}
        \Lambda \int_\Omega w^4 \dd x = 0.
    \end{equation}
    As $w \neq 0$, this gives $\Lambda = 0$, yielding $- \Delta w + \beta w - w^3 = 0$. Hence, $w \in \mathcal{N}^\prime$, and this completes the proof of $d = \inf_{\mathcal{N}} J$. By elliptic regularity, one has $w \in H^2(\Omega)$.
\end{proof}

We now prove Proposition \ref{prop:Q_nonconstant_lambda_1_smaller_than_2}. Hence, we assume that $\partial \Omega = \emptyset$. Let $0 = \lambda_0 < \lambda_1 < \cdots$ denote the eigenvalues of $-\Delta$, and let $(\varphi_k)_{k \geq 0}$ be the corresponding $L^2(\Omega)$-orthonormal basis of eigenfunctions, with, for instance, $\varphi_0 = \frac{1}{|\Omega|^{1/2}}$. We assume that $\lambda_1 < 2$, and we prove that $d < J_\KG(1)$.

\paragraph{Step 1: Local estimation of $J_\KG$ near $1$.}
Consider $h \in H^2(\Omega)$. A straightforward calculation gives
\begin{align}
    J_\KG(1 + h) \ & = \frac{1}{2} \int_\Omega \left\vert \nabla h \right\vert^2 \dd x + \frac{1}{2} \int_\Omega \left( 1 +  h \right)^2 \dd x - \frac{1}{4} \int_\Omega \left( 1 + h \right)^4 \dd x \\
    \ & = J_\KG(1) + \frac{1}{2} \int_\Omega \left\vert \nabla h \right\vert^2 \dd x - \int_\Omega h^2 \dd x - \int_\Omega h^3 \dd x - \frac{1}{4} \int_\Omega h^4 \dd x \\
    & = J_\KG(1) + \frac{1}{2} \left\langle Lh, h \right\rangle_{L^2(\Omega)}  - \int_\Omega h^3 \dd x - \frac{1}{4} \int_\Omega h^4 \dd x, \label{eq:proof:prop:Q_nonconstant_1}
\end{align}
where $Lh := - \Delta h - 2h$. The idea of the proof is to choose a small $h$ in the linear span of $\varphi_0$ and $\varphi_1$ such that 
\begin{equation}
    J_\KG(1 + h) \approx J_\KG(1) + \frac{1}{2} \left\langle Lh, h \right\rangle_{L^2(\Omega)} < J_\KG(1),
\end{equation}
and such that $1+h \in \mathcal{N}$, that is, $K(1+h) = 0$. This will show that the stationary solutions $\pm 1$ are not minimizers of the energy on the Nehari manifold $\mathcal{N}$, and thus that the energy of ground states is strictly below $J_\KG(1)$.

Consider $(\alpha, \beta) \in [-1, 1]^2$, and set $h_{\alpha, \beta} = \alpha \varphi_0 + \beta \varphi_1$. On the one hand, one has
\begin{equation}\label{eq:proof:prop:Q_nonconstant_2}
    \frac{1}{2} \left\langle Lh_{\alpha, \beta}, h_{\alpha, \beta} \right\rangle_{L^2(\Omega)} = \frac{1}{2} \int_\Omega \left( \beta \lambda_1 \varphi_1 - 2 \left( \alpha \varphi_0 + \beta \varphi_1 \right) \right) h_{\alpha, \beta} \dd x = - \alpha^2 + (\lambda_1 - 2) \frac{\beta^2}{2},
\end{equation}
and since $\lambda_1 < 2$, it implies that $\left\langle Lh_{\alpha, \beta}, h_{\alpha, \beta} \right\rangle < 0$ if $(\alpha, \beta) \neq (0,0)$. On the other hand, the remainder term that appears in \eqref{eq:proof:prop:Q_nonconstant_1} satisfies
\begin{equation}
    \left\vert \int_\Omega h_{\alpha, \beta}^3 \dd x + \frac{1}{4} \int_\Omega h_{\alpha, \beta}^4 \dd x \right\vert \leq C_1 \left( \vert \alpha \vert^3 + \vert \beta \vert^3 \right),
\end{equation}
for some $C_1 > 0$ which depends only on $\Omega$. Together with \eqref{eq:proof:prop:Q_nonconstant_1} and \eqref{eq:proof:prop:Q_nonconstant_2}, it gives
\begin{equation}\label{eq:proof:prop:Q_nonconstant_3}
    J_\KG(1+h_{\alpha, \beta}) \leq J_\KG(1) - \alpha^2 + (\lambda_1 - 2) \frac{\beta^2}{2} + C_1 \left( \vert \alpha \vert^3 + \vert \beta \vert^3 \right),
\end{equation}
and so
\begin{equation}\label{eq:proof:prop:Q_nonconstant_4}
    J_\KG(1+h_{\alpha, \beta}) < J_\KG(1) \quad \text{ if } 0 \leq \vert \alpha \vert < \frac{1}{C_1} \text{ and } 0 < \vert \beta \vert < \frac{2 - \lambda_1}{2 C_1}.
\end{equation}

\paragraph{Step 2: Choice of $h$ such that $1+h \in \mathcal{N}$.} One has
\begin{align}
    K(1 + h) \ & = \int_\Omega \left\vert \nabla h \right\vert^2 \dd x + \int_\Omega \left( 1 +  h \right)^2 \dd x - \int_\Omega \left( 1 + h \right)^4 \dd x \\
    & = - \int_\Omega  (\Delta h) h \dd x - 2 \int_\Omega h \dd x - 5 \int_\Omega h^2 \dd x - 4 \int_\Omega h^3 \dd x - \int_\Omega h^4 \dd x,
\end{align}
yielding
\begin{align}
    K(1 + h_{\alpha, \beta}) = - 2 \vert \Omega \vert^{\frac{1}{2}} \alpha + (\lambda_1 - 5) \beta^2 - 5 \alpha^2- 4 \int_\Omega h_{\alpha, \beta}^3 \dd x - \int_\Omega h_{\alpha, \beta}^4 \dd x. \label{eq:proof:prop:Q_nonconstant_5}
\end{align}
Set $a = 2 \vert \Omega \vert^{\frac{1}{2}} > 0$, $b = 5 - \lambda_1 > 0$, and 
\begin{equation}
    R(h_{\alpha, \beta}) = 5 \alpha^2 + 4 \int_\Omega h_{\alpha, \beta}^3 \dd x + \int_\Omega h_{\alpha, \beta}^4 \dd x,
\end{equation}
so that $K(1+h_{\alpha, \beta}) = -a \alpha - b \beta^2 - R(h_{\alpha, \beta})$. Note that $\vert R(h_{\alpha, \beta}) \vert \leq C_2 \left( \alpha^2 + \beta^3 \right)$, for some $C_2 > 0$ which depends only on $\Omega$. If $\alpha = 0$ and $0 \leq \beta \leq \frac{b}{C_2}$, then
\begin{equation}
    K(1+h_{\alpha, \beta}) \leq -b \beta^2 + C_2 \beta^3 = -b \beta^2 \left(1 - \frac{C_2 \beta}{b} \right) \leq 0.
\end{equation}
Similarly, there exists a small constant $c_3 > 0$, which depends only on $\Omega$, such that if $\alpha = - \frac{2 b \beta^2}{a}$ and $0 \leq \beta \leq c_3$, then
\begin{equation}
    K(1+h_{\alpha, \beta}) \geq b \beta^2 - C_2 \left( \beta^3 + \frac{4 b^2 \beta^4}{a^2} \right) \geq 0.
\end{equation}
Set $c_4 := \min \left(\frac{b}{C_2}, c_3 \right)$. By the intermediate value theorem, for all $\beta \in [0, c_4]$, there exists $\alpha = \alpha(\beta) \in \left[-\frac{2 b \beta^2}{a},0\right]$ such that $K(1+h_{\alpha, \beta}) = 0$. If, in addition, $\beta \neq 0$ and $\beta$ is sufficiently small, then $J_\KG(1+h_{\alpha, \beta}) < J_\KG(1)$ by \eqref{eq:proof:prop:Q_nonconstant_4}. This completes the proof of Proposition \ref{prop:Q_nonconstant_lambda_1_smaller_than_2}.


\subsection{Construction of solutions of the Klein-Gordon equation with interesting behaviours}

Here, we prove Theorem \ref{thm:construction_solution_particuliere}. Assume that $\partial \Omega = \emptyset$ and that $\gamma$ is a positive constant. Let $U_1 = U_1(\gamma) > 0$ be given by Lemma \ref{lem:gamma_fix_u_1_varies}. Note that the solution $u$ of \eqref{duffing} with initial data $\left( -1, U_1 \right)$ and damping $\gamma$ is a space-independent solution of \eqref{KG}, with constant initial data and constant damping, which satisfies 
\begin{equation}
    E_\KG \left(u(0), \partial_t u(0) \right) = \vert \Omega \vert E\left(u(0), \partial_t u(0) \right) \geq \frac{\vert \Omega \vert}{4} = \vert \Omega \vert E(1,0) = J_\KG(1),
\end{equation}
and it converges to $(1, 0)$ in $H_0^1(\Omega) \times L^2(\Omega)$ at an exponential rate. The same is true with $(1, 0)$ replaced by $(-1, 0)$ for the solution with initial data $\left( 1, -U_1 \right)$. It remains to prove the existence of the sets $\mathcal{O}_1$ and $\mathcal{O}_2$.

\paragraph{Construction of the set $\mathcal{O}_2$.} Let $u$ be the solution of \eqref{duffing} with initial data $\left( -1, \frac{U_1}{2} \right)$ and damping $\gamma$. By Lemma \ref{lem:gamma_fix_u_1_varies}, $u$ satisfies \eqref{eq:lem:gamma_fix_u_1_varies_1}. We also write $u$ for the solution $u(t,x) = u(t)$ of \eqref{KG} with constant initial data $x \mapsto \left( -1, \frac{U_1}{2} \right)$ and constant damping $\gamma$.

Since $E_\KG$ and $K$ (defined by \eqref{eq_def_K}) are continuous, there exists $\delta_0$ such that for all $(v_0, v_1) \in H_0^1(\Omega) \times L^2(\Omega)$, one has
\begin{equation}
    \left\Vert (v_0, v_1) \right\Vert_{H_0^1(\Omega) \times L^2(\Omega)} \leq \delta_0 \quad \Longrightarrow \quad (v_0, v_1) \in \mathcal{K}_\KG^+.
\end{equation}
Let $T$ be such that 
\begin{equation}
    \left\Vert \left( u(T), \partial_t u(T) \right) \right\Vert_{H_0^1(\Omega) \times L^2(\Omega)} \leq \frac{\delta_0}{2}. 
\end{equation}
By Lemma \ref{lem:source_to_sol}, there exists a neighbourhood $\mathcal{O}$ of the constant function $\left( -1, \frac{U_1}{2} \right)$ in $H_0^1(\Omega) \times L^2(\Omega)$, and there exists $C > 0$, which depends on $u$ and $T$, such that if $(v_0, v_1) \in \mathcal{O}$, then the solution $v$ of \eqref{KG} with initial data $\left( v_0, v_1 \right)$ and damping $\gamma$ satisfies
\begin{equation}
    \sup_{t \in [0, T]} \left\Vert \left( u, \partial_t u \right)(t) -  \left( v, \partial_t v \right)(t) \right\Vert_{H_0^1(\Omega) \times L^2(\Omega)} \leq C \left\Vert \left( -1, \frac{U_1}{2} \right) -  \left( v_0, v_1 \right) \right\Vert_{H_0^1(\Omega) \times L^2(\Omega)} .
\end{equation}
In particular, if in addition
\begin{equation}
    \left\Vert \left( -1, \frac{U_1}{2} \right) -  \left( v_0, v_1 \right) \right\Vert_{H_0^1(\Omega) \times L^2(\Omega)} \leq \frac{\delta_0}{2C},
\end{equation}
then $\left( v(T), \partial_t v(T) \right) \in \mathcal{K}_\KG^+$, and hence $\left( v, \partial_t v \right)$ tends to zero at an exponential rate by Theorem 4 of \cite{perrin2024damped}. 

\paragraph{Construction of the set $\mathcal{O}_1$.} Let $u$ be the solution of \eqref{duffing} with initial data $\left( -1, 2 U_1 \right)$ and damping $\gamma$. As above, we also write $u$ for the corresponding solution of \eqref{KG}. By Lemma \ref{lem:gamma_fix_u_1_varies}, $u$ blows up in finite time. Hence, by Theorem \ref{thm_duffing_N}, there exists $T > 0$ such that $u$ is defined on $[0, T]$ and $E_u(T) < 0$. Arguing as in the construction of $\mathcal{O}_2$, one proves that there exists an open set $\mathcal{O}_1$ such that for all $(v_0, v_1) \in \mathcal{O}_1$, the solution $v$ of \eqref{KG} with initial data $\left( v_0, v_1 \right)$ and damping $\gamma$ is defined on $[0, T]$ and satisfies 
\begin{equation}
    E_\KG\left( v(T), \partial_t v(T) \right) < 0.
\end{equation}
It implies $K(v(T)) < 0$, and hence $\left( v(T), \partial_t v(T) \right) \in \mathcal{K}_\KG^-$. By Theorem 2 of \cite{perrin2024damped}, $v$ blows up in finite time. This completes the proof of Theorem~\ref{thm:construction_solution_particuliere}.

\section*{Acknowledgments}

The author acknowledges support from the Fondation Simone et Cino Del Duca -- Institut de France, the Fondation Université de Rennes and from grant ANR-11-LABX-0020 (Labex Lebesgue).

\printbibliography

@article{bahri_existence_1997,
	title        = {On the existence of a positive solution of semilinear elliptic equations in unbounded domains},
	author       = {Bahri, A. and Lions, P.-L.},
	volume       = {14},
	number       = {3},
	pages        = {365--413},
	doi          = {10.1016/s0294-1449(97)80142-4},
	issn         = {0294-1449},
	url          = {https://ems.press/journals/aihpc/articles/4398274},
	journaltitle = {Annales de l'Institut Henri Poincaré C},
	date         = {1997-01-01},
	langid       = {english}
}

@article{Bandle1989,
	title        = {Semilinear elliptic problems in annular domains},
	author       = {Bandle,  C. and Kwong,  M. K.},
	year         = {1989},
	month        = mar,
	journal      = {ZAMP Zeitschrift für angewandte Mathematik und Physik},
	publisher    = {Springer Science and Business Media LLC},
	volume       = {40},
	number       = {2},
	pages        = {245–257},
	doi          = {10.1007/bf00945001},
	issn         = {1420-9039},
	url          = {http://dx.doi.org/10.1007/BF00945001}
}

@article{BLR,
	title        = {Sharp sufficient conditions for the observation, control, and stabilization of waves from the boundary},
	author       = {Bardos, C. and Lebeau, G. and Rauch, J.},
	year         = {1992},
	journal      = {SIAM J. Control Optim.},
	volume       = {30},
	number       = {5},
	pages        = {1024--1065},
	doi          = {10.1137/0330055},
	fjournal     = {SIAM Journal on Control and Optimization}
}

@article{Benci1987,
	title        = {Positive solutions of some nonlinear elliptic problems in exterior domains},
	author       = {Benci,  V. and Cerami,  G.},
	year         = {1987},
	journal      = {Archive for Rational Mechanics and Analysis},
	publisher    = {Springer Science and Business Media LLC},
	volume       = {99},
	number       = {4},
	pages        = {283–300},
	doi          = {10.1007/bf00282048},
	issn         = {1432-0673},
	url          = {http://dx.doi.org/10.1007/BF00282048}
}

@article{berestycki_methode_1980,
	title        = {Une methode locale pour l’existence de solutions positives de problemes semi-lineaires elliptiques dans $R^N$},
	author       = {Berestycki, H. and Lions, P.-L.},
	volume       = {38},
	number       = {1},
	pages        = {144--187},
	doi          = {10.1007/BF03033880},
	issn         = {1565-8538},
	url          = {https://doi.org/10.1007/BF03033880},
	journaltitle = {Journal d’Analyse Mathématique},
	shortjournal = {J. Anal. Math.},
	date         = {1980-12-01},
	langid       = {french}
}

@article{brezis_positive_1983,
	title        = {Positive solutions of nonlinear elliptic equations involving critical sobolev exponents},
	author       = {Brezis, H. and Nirenberg, L.},
	year         = {1983},
	journal      = {Communications on Pure and Applied Mathematics},
	volume       = {36},
	number       = {4},
	pages        = {437--477},
	doi          = {10.1002/cpa.3160360405},
	issn         = {1097-0312},
	url          = {https://onlinelibrary.wiley.com/doi/abs/10.1002/cpa.3160360405},
	language     = {en}
}

@article{BurqRaugelSchlag,
	title        = {Long time dynamics for damped Klein-Gordon equations},
	author       = {Burq, N. and Raugel, G. and Schlag, W.},
	year         = {2017},
	journal      = {Annales scientifiques de l{\textquotesingle}{\'{E}}cole normale sup{\'{e}}rieure},
	publisher    = {Societe Mathematique de France},
	volume       = {50},
	number       = {6},
	pages        = {1447--1498},
	doi          = {10.24033/asens.2349}
}

@article{byeon_addendum_2001,
	title        = {{ADDENDUM}: Volume 163, Number 2 (2000), in the article “Effect of Symmetry to the Structure of Positive Solutions in Nonlinear Eliptic Problems,” by Jaeyoung Byeon, pages 429–474 ()},
	shorttitle   = {{ADDENDUM}},
	author       = {Byeon, J.},
	volume       = {172},
	number       = {2},
	pages        = {445--447},
	doi          = {10.1006/jdeq.2001.4040},
	issn         = {0022-0396},
	url          = {https://www.sciencedirect.com/science/article/pii/S0022039601940409},
	journaltitle = {Journal of Differential Equations},
	shortjournal = {Journal of Differential Equations},
	date         = {2001-05-20}
}

@article{byeon_effect_2000,
	title        = {Effect of Symmetry to the Structure of Positive Solutions in Nonlinear Eliptic Problems},
	author       = {Byeon, J.},
	volume       = {163},
	number       = {2},
	pages        = {429--474},
	doi          = {10.1006/jdeq.1999.3737},
	issn         = {0022-0396},
	url          = {https://www.sciencedirect.com/science/article/pii/S0022039699937373},
	journaltitle = {Journal of Differential Equations},
	shortjournal = {Journal of Differential Equations},
	date         = {2000-05-20}
}

@article{byeon_existence_1997,
	title        = {Existence of Many Nonequivalent Nonradial Positive Solutions of Semilinear Elliptic Equations on Three-Dimensional Annuli},
	author       = {Byeon, J.},
	volume       = {136},
	number       = {1},
	pages        = {136--165},
	doi          = {10.1006/jdeq.1996.3241},
	issn         = {0022-0396},
	url          = {https://www.sciencedirect.com/science/article/pii/S0022039696932416},
	journaltitle = {Journal of Differential Equations},
	shortjournal = {Journal of Differential Equations},
	date         = {1997-05-01}
}

@article{chen_uniqueness_1991,
	title        = {Uniqueness of the ground state solutions of $\Delta u+f(u)=0$ in $\mathbb{R}^n$, $n\geq3$},
	author       = {Chen, C.-C. and Lin, C.-S.},
	volume       = {16},
	number       = {8},
	pages        = {1549--1572},
	doi          = {10.1080/03605309108820811},
	issn         = {0360-5302},
	url          = {https://doi.org/10.1080/03605309108820811},
	note         = {Publisher: Taylor \& Francis \_eprint: https://doi.org/10.1080/03605309108820811},
	journaltitle = {Communications in Partial Differential Equations},
	date         = {1991-01-01}
}

@article{Coffman1972,
	title        = {Uniqueness of the ground state solution for $\Delta u-u+u^3=0$ and a variational characterization of other solutions},
	author       = {Coffman,  C. V.},
	year         = {1972},
	month        = jan,
	journal      = {Archive for Rational Mechanics and Analysis},
	publisher    = {Springer Science and Business Media LLC},
	volume       = {46},
	number       = {2},
	pages        = {81–95},
	doi          = {10.1007/bf00250684},
	issn         = {1432-0673},
	url          = {http://dx.doi.org/10.1007/BF00250684}
}

@article{Duyckaerts2024,
	title        = {Global Solutions with Asymptotic Self-Similar Behaviour for the Cubic Wave Equation},
	author       = {Duyckaerts,  T. and Negro,  G.},
	year         = {2024},
	month        = mar,
	journal      = {Communications in Mathematical Physics},
	publisher    = {Springer Science and Business Media LLC},
	volume       = {405},
	number       = {3},
	doi          = {10.1007/s00220-024-04962-3},
	issn         = {1432-0916},
	url          = {http://dx.doi.org/10.1007/s00220-024-04962-3}
}

@article{DuyckaertsMerle,
	title        = {{Dynamics of Threshold Solutions for Energy-Critical Wave Equation}},
	author       = {Duyckaerts, T. and Merle, F.},
	year         = {2008},
	month        = {01},
	journal      = {International Mathematics Research Papers},
	volume       = {2008},
	pages        = {rpn002},
	doi          = {10.1093/imrp/rpn002},
	issn         = {1687-3017}
}

@article{Esteban1982,
	title        = {Existence and non-existence results for semilinear elliptic problems in unbounded domains},
	author       = {Esteban,  M. J. and Lions,  P.-L.},
	year         = {1982},
	journal      = {Proceedings of the Royal Society of Edinburgh: Section A Mathematics},
	publisher    = {Cambridge University Press (CUP)},
	volume       = {93},
	number       = {1–2},
	pages        = {1–14},
	doi          = {10.1017/s0308210500031607},
	issn         = {1473-7124},
	url          = {http://dx.doi.org/10.1017/S0308210500031607}
}

@article{Duffing_defocusing,
	title        = {Monotonous property of non-oscillations of the damped {Duffing}'s equation},
	author       = {Feng, Z.},
	year         = {2006},
	journal      = {Chaos Solitons Fractals},
	volume       = {28},
	number       = {2},
	pages        = {463--471},
	doi          = {10.1016/j.chaos.2005.07.006},
	issn         = {0960-0779},
	fjournal     = {Chaos, Solitons and Fractals},
	language     = {English},
	zbmath       = {5005529},
	zbl          = {1097.34026}
}

@article{gazzola_global_2006,
	title        = {Global solutions and finite time blow up for damped semilinear wave equations},
	author       = {Gazzola, F. and Squassina, M.},
	year         = {2006},
	month        = apr,
	journal      = {Annales de l'Institut Henri Poincaré C},
	volume       = {23},
	number       = {2},
	pages        = {185--207},
	doi          = {10.1016/j.anihpc.2005.02.007},
	issn         = {0294-1449},
	url          = {https://ems.press/journals/aihpc/articles/4077621},
	language     = {en}
}

@article{Gidas1979,
	title        = {Symmetry and related properties via the maximum principle},
	author       = {Gidas,  B. and Ni,  W.-M. and Nirenberg,  L.},
	year         = {1979},
	month        = oct,
	journal      = {Communications in Mathematical Physics},
	publisher    = {Springer Science and Business Media LLC},
	volume       = {68},
	number       = {3},
	pages        = {209–243},
	doi          = {10.1007/bf01221125},
	issn         = {1432-0916},
	url          = {http://dx.doi.org/10.1007/BF01221125}
}

@article{grillakis_existence_1990,
	title        = {Existence of nodal solutions of semilinear equations in $\mathbb{R}^N$},
	author       = {Grillakis, M.},
	volume       = {85},
	number       = {2},
	pages        = {367--400},
	doi          = {10.1016/0022-0396(90)90121-5},
	issn         = {0022-0396},
	url          = {https://www.sciencedirect.com/science/article/pii/0022039690901215},
	journaltitle = {Journal of Differential Equations},
	shortjournal = {Journal of Differential Equations},
	date         = {1990-06-01}
}

@article{IbrahimMasmoudiNakanishi,
	title        = {Scattering threshold for the focusing nonlinear Klein{\textendash}Gordon equation},
	author       = {Ibrahim, S. and Masmoudi, N. and Nakanishi, K.},
	year         = {2011},
	month        = dec,
	journal      = {Analysis \& PDE},
	publisher    = {Mathematical Sciences Publishers},
	volume       = {4},
	number       = {3},
	pages        = {405--460},
	doi          = {10.2140/apde.2011.4.405},
	url          = {https://doi.org/10.2140/apde.2011.4.405}
}

@article{johnson_singular_1994,
	title        = {Singular solutions of the elliptic equation $\delta u - u + u^p=0$},
	author       = {Johnson, R. A. and Pan, X. and Yi, Y.},
	volume       = {166},
	number       = {1},
	pages        = {203--225},
	doi          = {10.1007/BF01765635},
	issn         = {1618-1891},
	url          = {https://doi.org/10.1007/BF01765635},
	journaltitle = {Annali di Matematica Pura ed Applicata},
	shortjournal = {Annali di Matematica pura ed applicata},
	date         = {1994-12-01},
	langid       = {english}
}

@article{Joly-Laurent,
	title        = {{Stabilization for the semilinear wave equation with geometric control condition}},
	author       = {Joly, R. and Laurent, C.},
	year         = {2013},
	journal      = {Analysis \& PDE},
	publisher    = {MSP},
	volume       = {6},
	number       = {5},
	pages        = {1089 -- 1119},
	doi          = {10.2140/apde.2013.6.1089}
}

@article{jones_infinitely_1986,
	title        = {On the Infinitely Many Solutions of a Semilinear Elliptic Equation},
	author       = {Jones, C. and Küpper, T.},
	volume       = {17},
	number       = {4},
	pages        = {803--835},
	doi          = {10.1137/0517059},
	issn         = {0036-1410},
	url          = {https://epubs.siam.org/doi/10.1137/0517059},
	note         = {Publisher: Society for Industrial and Applied Mathematics},
	journaltitle = {{SIAM} Journal on Mathematical Analysis},
	shortjournal = {{SIAM} J. Math. Anal.},
	date         = {1986-07}
}

@article{KenigMerle,
	title        = {Global well-posedness,  scattering and blow-up for the energy-critical focusing non-linear wave equation},
	author       = {Kenig, C. E. and Merle, F.},
	year         = {2008},
	journal      = {Acta Mathematica},
	publisher    = {International Press of Boston},
	volume       = {201},
	number       = {2},
	pages        = {147–212},
	doi          = {10.1007/s11511-008-0031-6},
	issn         = {0001-5962},
	url          = {http://dx.doi.org/10.1007/s11511-008-0031-6}
}

@article{krieger_global_2012-1D,
	title        = {Global dynamics above the ground state energy for the one-dimensional {NLKG} equation},
	author       = {Krieger, J. and Nakanishi, K. and Schlag, W.},
	year         = {2012},
	month        = oct,
	journal      = {Mathematische Zeitschrift},
	volume       = {272},
	number       = {1},
	pages        = {297--316},
	doi          = {10.1007/s00209-011-0934-3},
	issn         = {1432-1823},
	url          = {https://doi.org/10.1007/s00209-011-0934-3},
	language     = {en}
}

@article{krieger_global_2012,
	title        = {Global dynamics of the nonradial energy-critical wave equation above the ground state energy},
	author       = {Krieger, J. and Nakanishi, K. and Schlag, W.},
	year         = {2012},
	month        = dec,
	journal      = {Discrete and Continuous Dynamical Systems},
	volume       = {33},
	number       = {6},
	pages        = {2423--2450},
	doi          = {10.3934/dcds.2013.33.2423},
	issn         = {1078-0947},
	url          = {https://www.aimsciences.org/en/article/doi/10.3934/dcds.2013.33.2423},
	copyright    = {http://creativecommons.org/licenses/by/3.0/},
	note         = {Publisher: Discrete and Continuous Dynamical Systems},
	language     = {en}
}

@article{krieger_global_2013,
	title        = {Global dynamics away from the ground state for the energy-critical nonlinear wave equation},
	author       = {Krieger, J. and Nakanishi, K. and Schlag, W.},
	year         = {2013},
	journal      = {American Journal of Mathematics},
	volume       = {135},
	number       = {4},
	pages        = {935--965},
	issn         = {1080-6377},
	url          = {https://muse.jhu.edu/pub/1/article/514419},
	note         = {Publisher: Johns Hopkins University Press}
}

@article{krieger_threshold_2014,
	title        = {Threshold {Phenomenon} for the {Quintic} {Wave} {Equation} in {Three} {Dimensions}},
	author       = {Krieger, J. and Nakanishi, K. and Schlag, W.},
	year         = {2014},
	month        = apr,
	journal      = {Communications in Mathematical Physics},
	volume       = {327},
	number       = {1},
	pages        = {309--332},
	doi          = {10.1007/s00220-014-1900-9},
	issn         = {1432-0916},
	url          = {https://doi.org/10.1007/s00220-014-1900-9},
	language     = {en}
}

@misc{krieger_center-stable_2013,
	title        = {Center-stable manifold of the ground state in the energy space for the critical wave equation},
	author       = {Krieger, J. and Nakanishi, K. and Schlag, W.},
	year         = {2015},
	journal      = {Math. Ann.},
	volume       = {361},
	number       = {1-2},
	pages        = {1--50},
	doi          = {10.1007/s00208-014-1059-x},
	issn         = {0025-5831},
	url          = {hdl.handle.net/2066/141321},
	fjournal     = {Mathematische Annalen},
	language     = {English}
}

@article{Kwong1989,
	title        = {Uniqueness of positive solutions of $\Delta u-u+u^p=0$ in $R^n$},
	author       = {Kwong,  M. K.},
	year         = {1989},
	month        = sep,
	journal      = {Archive for Rational Mechanics and Analysis},
	publisher    = {Springer Science and Business Media LLC},
	volume       = {105},
	number       = {3},
	pages        = {243–266},
	doi          = {10.1007/bf00251502},
	issn         = {1432-0673},
	url          = {http://dx.doi.org/10.1007/BF00251502}
}

@article{LiZaho2020,
	title        = {Asymptotic Decomposition for Nonlinear Damped Klein-Gordon Equations},
	author       = {Li,  Z. and Zhao,  L.},
	year         = {2020},
	month        = may,
	journal      = {Journal of Mathematical Study},
	publisher    = {Global Science Press},
	volume       = {53},
	number       = {3},
	pages        = {329–352},
	doi          = {10.4208/jms.v53n3.20.06},
	issn         = {2096-9856},
	url          = {http://dx.doi.org/10.4208/jms.v53n3.20.06}
}

@article{li_monotonicity_1991,
	title        = {Monotonicity and Symmetry of Solutions of Fully Nonlinear Elliptic Equations on Bounded domains},
	author       = {Li, C.},
	volume       = {16},
	number       = {2},
	pages        = {491--526},
	doi          = {10.1080/03605309108820766},
	issn         = {0360-5302},
	url          = {https://doi.org/10.1080/03605309108820766},
	note         = {Publisher: Taylor \& Francis \_eprint: https://doi.org/10.1080/03605309108820766},
	journaltitle = {Communications in Partial Differential Equations},
	date         = {1991-01-01}
}

@article{lin_asymptotic_1995,
	title        = {Asymptotic Behavior of Positive Solutions to Semilinear Elliptic Equations on Expanding Annuli},
	author       = {Lin, S. S.},
	volume       = {120},
	number       = {2},
	pages        = {255--288},
	doi          = {10.1006/jdeq.1995.1112},
	issn         = {0022-0396},
	url          = {https://www.sciencedirect.com/science/article/pii/S0022039685711126},
	journaltitle = {Journal of Differential Equations},
	shortjournal = {Journal of Differential Equations},
	date         = {1995-08-01}
}

@article{lin_existence_1992,
	title        = {Existence of positive nonradial solutions for nonlinear elliptic equations in annular domains},
	author       = {Lin, S. S.},
	volume       = {332},
	number       = {2},
	pages        = {775--791},
	doi          = {10.1090/S0002-9947-1992-1055571-1},
	issn         = {0002-9947, 1088-6850},
	url          = {https://www.ams.org/tran/1992-332-02/S0002-9947-1992-1055571-1/},
	journaltitle = {Transactions of the American Mathematical Society},
	shortjournal = {Trans. Amer. Math. Soc.},
	date         = {1992},
	langid       = {english}
}

@article{mcleod_radial_1990,
	title        = {Radial solutions of $\delta u+ f (u)= 0$ with prescribed numbers of zeros},
	author       = {{McLeod}, K. and Troy, W. C. and Weissler, F. B.},
	volume       = {83},
	number       = {2},
	pages        = {368--378},
	url          = {https://www.sciencedirect.com/science/article/pii/002203969090063U},
	note         = {Publisher: Elsevier},
	journaltitle = {Journal of Differential Equations},
	date         = {1990}
}

@article{MelroseSjostrand,
	title        = {Singularities of boundary value problems. I},
	author       = {Melrose, R. B. and Sj\"{o}strand, J.},
	year         = {1978},
	month        = sep,
	journal      = {Communications on Pure and Applied Mathematics},
	publisher    = {Wiley},
	volume       = {31},
	number       = {5},
	pages        = {593--617},
	doi          = {10.1002/cpa.3160310504}
}

@article{NakanishiSchlagArticleRAD,
	title        = {Global dynamics above the ground state energy for the focusing nonlinear Klein–Gordon equation},
	author       = {Nakanishi, K. and Schlag, W.},
	year         = {2011},
	journal      = {Journal of Differential Equations},
	volume       = {250},
	number       = {5},
	pages        = {2299--2333},
	doi          = {https://doi.org/10.1016/j.jde.2010.10.027},
	issn         = {0022-0396},
	url          = {https://www.sciencedirect.com/science/article/pii/S0022039610004213}
}

@book{NakanishiSchlagBOOK,
	title        = {Invariant manifolds and dispersive {Hamiltonian} evolution equations},
	author       = {Nakanishi, K. and Schlag, W.},
	year         = {2011},
	publisher    = {Z{\"u}rich: European Mathematical Society (EMS)},
	series       = {Zur. Lect. Adv. Math.},
	doi          = {10.4171/095},
	isbn         = {978-3-03719-095-1},
	fseries      = {Zurich Lectures in Advanced Mathematics},
	language     = {English}
}

@article{NakanishiSchlagArticleNONRAD,
	title        = {Global dynamics above the ground state for the nonlinear {Klein--Gordon} equation without a radial assumption},
	author       = {Nakanishi, K. and Schlag, W.},
	year         = {2012},
	journal      = {Arch. Ration. Mech. Anal.},
	publisher    = {Springer Science and Business Media LLC},
	volume       = {203},
	number       = {3},
	pages        = {809--851},
	doi          = {10.1007/s00205-011-0462-7}
}

@article{Ni1983,
	title        = {Uniqueness of solutions of nonlinear Dirichlet problems},
	author       = {Ni,  W.-M.},
	year         = {1983},
	month        = nov,
	journal      = {Journal of Differential Equations},
	publisher    = {Elsevier BV},
	volume       = {50},
	number       = {2},
	pages        = {289–304},
	doi          = {10.1016/0022-0396(83)90079-7},
	issn         = {0022-0396},
	url          = {http://dx.doi.org/10.1016/0022-0396(83)90079-7}
}

@article{noussair_existence_1979,
	title        = {On the existence of solutions of nonlinear elliptic boundary value problems},
	author       = {Noussair, E. S.},
	volume       = {34},
	number       = {3},
	pages        = {482--495},
	doi          = {10.1016/0022-0396(79)90032-9},
	issn         = {0022-0396},
	url          = {https://www.sciencedirect.com/science/article/pii/0022039679900329},
	journaltitle = {Journal of Differential Equations},
	shortjournal = {Journal of Differential Equations},
	date         = {1979-12-01}
}

@article{Payne-Sattinger,
	title        = {Saddle points and instability of nonlinear hyperbolic equations},
	author       = {Payne, L. E. and Sattinger, D. H.},
	year         = {1975},
	journal      = {Israel Journal of Mathematics},
	volume       = {22},
	pages        = {273--303},
	doi          = {10.1007/BF02761595}
}

@misc{perrin2024damped,
	title        = {The damped focusing cubic wave equation on a bounded domain},
	author       = {Perrin, T.},
	year         = {2024},
	url          = {https://arxiv.org/abs/2310.12644},
	eprint       = {2310.12644},
	archiveprefix = {arXiv},
	primaryclass = {math.AP}
}

@article{pucci_global_1998,
	title        = {Global {Nonexistence} for {Abstract} {Evolution} {Equations} with {Positive} {Initial} {Energy}},
	author       = {Pucci, P. and Serrin, J.},
	year         = {1998},
	month        = nov,
	journal      = {Journal of Differential Equations},
	volume       = {150},
	number       = {1},
	pages        = {203--214},
	doi          = {10.1006/jdeq.1998.3477},
	issn         = {0022-0396},
	url          = {https://www.sciencedirect.com/science/article/pii/S0022039698934775}
}

@article{seda_properties_1990,
	title        = {On some properties of the solution of the differential equation $u''+ \frac{2u'}{r}=u-u^3$},
	author       = {Šeda, V. and Pekár, J.},
	volume       = {35},
	number       = {4},
	pages        = {315--336},
	issn         = {0373-6725},
	url          = {https://dml.cz/handle/10338.dmlcz/104413},
	note         = {Publisher: Institute of Mathematics, Academy of Sciences of the Czech Republic},
	journaltitle = {Aplikace matematiky},
	date         = {1990}
}

@article{tanaka_uniqueness_2016,
	title        = {Uniqueness of sign-changing radial solutions for $\delta u -u+ \vert u \vert^{p-1} u=0$ in some ball and annulus},
	author       = {Tanaka, S.},
	volume       = {439},
	number       = {1},
	pages        = {154--170},
	doi          = {10.1016/j.jmaa.2016.02.036},
	issn         = {0022-247X},
	url          = {https://www.sciencedirect.com/science/article/pii/S0022247X16001700},
	journaltitle = {Journal of Mathematical Analysis and Applications},
	shortjournal = {Journal of Mathematical Analysis and Applications},
	date         = {2016-07-01}
}

@article{tang_uniqueness_2003,
	title        = {Uniqueness of positive radial solutions for $\delta u -u+u^p=0$ on an annulus},
	author       = {Tang, M.},
	volume       = {189},
	number       = {1},
	pages        = {148--160},
	doi          = {10.1016/S0022-0396(02)00142-0},
	issn         = {0022-0396},
	url          = {https://www.sciencedirect.com/science/article/pii/S0022039602001420},
	journaltitle = {Journal of Differential Equations},
	shortjournal = {Journal of Differential Equations},
	date         = {2003-03-20}
}

@article{troy_existence_2005,
	title        = {The existence and uniqueness of bound-state solutions of a semi-linear equation},
	author       = {Troy, W. C.},
	doi          = {10.1098/rspa.2005.1482},
	url          = {https://royalsocietypublishing.org/doi/10.1098/rspa.2005.1482},
	note         = {Publisher: The Royal {SocietyLondon}},
	rights       = {© 2005 The Royal Society},
	journaltitle = {Proceedings of the Royal Society A: Mathematical, Physical and Engineering Sciences},
	date         = {2005-08-04}
}

@article{vitillaro_global_1999,
	title        = {Global {Nonexistence} {Theorems} for a {Class} of {Evolution} {Equations} with {Dissipation}},
	author       = {Vitillaro, E.},
	year         = {1999},
	month        = oct,
	journal      = {Archive for Rational Mechanics and Analysis},
	volume       = {149},
	number       = {2},
	pages        = {155--182},
	doi          = {10.1007/s002050050171},
	issn         = {1432-0673},
	url          = {https://doi.org/10.1007/s002050050171},
	language     = {en}
}

@article{yarur_uniqueness_2009,
	title        = {On the uniqueness of the second bound state solution of a semilinear equation},
	author       = {Yarur, C. S. and Cortázar, C. and García-Huidobro, M.},
	volume       = {26},
	number       = {6},
	pages        = {2091--2110},
	doi          = {10.1016/j.anihpc.2009.01.004},
	issn         = {0294-1449},
	url          = {https://ems.press/journals/aihpc/articles/4077942},
	journaltitle = {Annales de l'Institut Henri Poincaré C},
	date         = {2009-12-01},
	langid       = {english}
}

@article{zou_effect_1994,
	title        = {On the effect of the domain geometry on uniqueness of positive solutions of $\Delta u + u^p = 0$},
	author       = {Zou, H.},
	volume       = {21},
	number       = {3},
	pages        = {343--356},
	issn         = {2036-2145},
	url          = {http://www.numdam.org/item/?id=ASNSP_1994_4_21_3_343_0},
	journaltitle = {Annali della Scuola Normale Superiore di Pisa - Classe di Scienze},
	date         = {1994},
	langid       = {french}
}

\noindent
\textsc{Perrin Thomas:} \texttt{thomas.perrin@ens-rennes.fr}

\end{document}